\documentclass[10.5pt]{article}
\usepackage[leqno]{amsmath}
\usepackage{amsfonts}
\usepackage{graphicx}

\usepackage{amsmath}
\usepackage{amssymb}
\usepackage{latexsym}
\usepackage{amsmath, amsfonts,amssymb, amsthm, enumerate, euscript,makeidx,color,mathrsfs}
\usepackage{tikz}
\usetikzlibrary{calc, shadings, shadows, shapes.arrows}

\def\sqr#1#2{{\vcenter{\vbox{\hrule height.#2pt
              \hbox{\vrule width.#2pt height#1pt \kern#1pt \vrule width.#2pt}
              \hrule height.#2pt}}}}

\def\ol{\overline}

\def\3n{\negthinspace \negthinspace \negthinspace }
\def\2n{\negthinspace \negthinspace }
\def\1n{\negthinspace }

\def\dbR{\mathbb{R}}
\def\dbS{\mathbb{S}}

\def\sA{{\mathscr A}}

\def\sD{{\mathscr D}}

\def\sF{{\mathscr F}}

\def\sM{{\mathscr M}}

\def\sU{{\mathscr U}}
\def\sV{{\mathscr V}}

\def\sY{{\mathscr Y}}

\def\={\buildrel \triangle \over =}

\def\ds{\displaystyle}

\def\ns{\noalign{\ss}}
\def\a{\alpha}
\def\b{\beta}
\def\g{\gamma}
\def\d{\delta}
\def\e{\varepsilon}

\def\l{\lambda}
\def\m{\mu}

\def\si{\sigma}

\def\f{\varphi}
\def\th{\theta}
\def\o{\omega}

\def\G{\Gamma}
\def\D{\Delta}

\def\L{\Lambda}
\def\Si{\Sigma}

\def\O{\Omega}

\def\cH{{\cal H}}

\def\cL{{\cal L}}

\def\cQ{{\cal Q}}

\def\cl{{\cal l}}
\def\BF{{\bf F}}

\def\BP{{\bf P}}

\def\Bf{{\bf f}}

\def\no{\noindent}

\def\ss{\smallskip}
\def\ms{\medskip}

\def\q{\quad}
\def\qq{\qquad}
\def\hb{\hbox}

\def\liminf{\mathop{\underline{\rm lim}}}

\def\da{\mathop{\downarrow}}
\def\Ra{\mathop{\Rightarrow}}

\def\lan{\mathop{\langle}}
\def\ran{\mathop{\rangle}}

\def\hcon{\stackrel{H}{\longrightarrow}}

\def\pa{\partial}
\def\na{\nabla}
\def\h{\widehat}
\def\wt{\widetilde}

\def\cd{\cdot}
\def\cds{\cdots}

\def\ae{\hbox{\rm a.e.{ }}}

\def\sgn{\hbox{\rm sgn$\,$}}

\def\supp{\hbox{\rm supp$\,$}}

\def\tr{\hbox{\rm tr$\,$}}

\def\cl{\overline}
\def\co{\mathop{{\rm co}}}

\def\les{\leqslant}
\def\ges{\geqslant}

\def\({\Big (}
\def\){\Big )}
\def\[{\Big[}
\def\]{\Big]}
\def\bde{\begin{definition}\label}
\def\ede{\end{definition}}
\def\be{\begin{equation}}
\def\bel{\begin{equation}\label}
\def\ee{\end{equation}}
\def\bt{\begin{theorem}\label}
\def\et{\end{theorem}}
\def\bc{\begin{corollary}\label}
\def\ec{\end{corollary}}
\def\bl{\begin{lemma}\label}
\def\el{\end{lemma}}
\def\bp{\begin{proposition}\label}
\def\ep{\end{proposition}}
\def\bex{\begin{example}\label}
\def\ex{\end{example}}
\def\bas{\begin{assumption}}
\def\eas{\end{assumption}}
\def\br{\begin{remark}\label}
\def\er{\end{remark}}
\def\ba{\begin{array}}
\def\ea{\end{array}}
\def\ed{\end{document}}

\def\square#1{\vbox{\hrule\hbox{\vrule height#1%
     \kern#1\vrule}\hrule}}
\def\rectangle#1#2{\vbox{\hrule\hbox{\vrule height#1%
     \kern#2\vrule}\hrule}}

\font\tenbb=msbm10 \font\sevenbb=msbm7 \font\fivebb=msbm5

\newfam\bbfam
\scriptscriptfont\bbfam=\fivebb \textfont\bbfam=\tenbb
\scriptfont\bbfam=\sevenbb

\newcommand{\sti}[1]{\mbox{ strongly in}\, #1}
\newcommand{\wi}[1]{\mbox{ weakly  in}\, #1}

\def\endpf{\hfill$\Box$\vspace{0.4cm}}

\def\eqon{ \, {\rm on } \, \,}
\def\eqin{ \, {\rm in } \, \,}
\def\eqae{ \, {\rm a.e. } \,\, }

\newcommand{\set}[1]{\left\{#1\right\}}
\newcommand{\ip}[1]{\left\langle #1\right\rangle}

\newcommand{\ppmatrix}[1]{\begin{pmatrix}#1\end{pmatrix}}

\def\all{  \, \forall \, }
\def\RT{\mathrm{T}}

\def\uGs{{\underline \sigma}}
\def\ula{{\underline a}}
\def\uly{{\underline y}}

\newtheorem{definition}{Definition}[section]
\newtheorem{lemma}[definition]{Lemma}
\newtheorem{theorem}[definition]{Theorem}
\newtheorem{corollary}[definition]{Corollary}
\newtheorem{proposition}[definition]{Proposition}
\newtheorem{remark}[definition]{Remark}
\newtheorem{example}[definition]{Example}
\newtheorem{conjecture}[definition]{Conjecture}
\newtheorem{assumption}[definition]{Assumption}

\if{
\newtheorem{lemma}{Lemma}[section]
\newtheorem{remark}{Remark}[section]
\newtheorem{example}{Example}[section]
\newtheorem{theorem}{Theorem}[section]
\newtheorem{corollary}{Corollary}[section]

\newtheorem{definition}{Definition}[section]
\newtheorem{proposition}{Proposition}[section]
\newtheorem{assumption}{Assumption}[section]
}\fi

\makeatletter
   
   \@addtoreset{equation}{section}
\makeatother

\begin{document}

\title{\bf Optimization of the Principal Eigenvalue\\ for Elliptic Operators
\thanks{This work is supported in part by NSFC Grant 11771097 and by NSF Grant DMS-1812921.}}
\author{Hongwei Lou\thanks{School of Mathematical Sciences, Fudan University,
Shanghai 200433, China; Email: hwlou@fudan.edu.cn}\, \ \ and \ \
Jiongmin Yong
\thanks{Department of Mathematics, University of Central Florida, Orlando, FL 32816, USA;
Email:jiongmin.yong @ucf.edu.}
 }
\maketitle

\no\bf Abstract: \rm Maximization and minimization problems of the principle eigenvalue for divergence form second order elliptic operators with the Dirichlet boundary condition are considered. The principal eigen map of such elliptic operators is introduced and some basic properties of this map, including continuity, concavity, and differentiability with respect to the parameter in the diffusibility matrix, are established. For maximization problem, the admissible control set is convexified to get the existence of an optimal convexified relaxed solution. Whereas, for minimization problem, the relaxation of the problem under $H$-convergence is introduced to get an optimal $H$-relaxed solution for certain interesting special cases. Some necessary optimality conditions are presented for both problems and a couple of illustrative examples are presented as well.

\ms

\no\bf Keywords: \rm elliptic operator, principal eigenvalue,
normalized principal eigenfunction, principal eigen map, $H$-convergence,
necessary optimality condition.

\ms

\no\bf AMS Mathematics Subject Classification. \rm 35J15, 35P05,
47A75, 49K20, 49J20.

\section{Introduction}

Consider a heat conduct problem in a bounded domain $\O\subseteq\dbR^n$. Suppose $\O$ is occupied by a certain type of medium with (not necessarily isotropic) uniformly elliptic {\it diffusibility matrix} $a(\cd)\equiv\big(a_{ij}(\cd)\big)$. Let $y(t,x)$ be the temperature of the body at $(t,x)$. Then, in the case that there is neither source nor sink of the heat in the domain, and the temperature is set to be a fixed level (say, 0, for simplicity) at the boundary $\pa\O$, the (temperature) function $y(\cd\,,\cd)\equiv y(\cd\,,\cd\,;y_0(\cd))$ will be the weak solution to the following parabolic equation:
\bel{1.1}\left\{\2n\ba{ll}
\ds y_t(t,x)=\nabla\cd\big(a(x)\nabla
y(t,x)\big),\qq(t,x)\in[0,\infty)\times\O,\\
\ns\ds y(t,x)=0,\qq(t,x)\in[0,\infty)\times\pa\O,\\
\ns\ds y(0,x)=y_0(x),\qq x\in\O.\ea\right.\ee
Consequently,
\bel{1.2}\ba{ll}
\ds{d\over dt}\|y(t,\cd)\|_2^2=2\int_\O y(t,x)\nabla\cd\big(a(x)\nabla y(t,x)\big)dx
=-2\int_\O\lan a(x)\nabla y(t,x),\nabla y(t,x)\ran dx\\
[4mm]
\ds\qq\qq\q\;\les-2\l_{a(\cd)}\int_\O|y(t,x)|^2dx\equiv-2\l_{a(\cd)}\|y(t,\cd)\|_2^2,\ea\ee
where $\|\cd\|_2=\|\cd\|_{L^2(\O)}$ is the usual $L^2(\O)$-norm, and $\l_{a(\cd)}>0$ is the {\it smallest
eigenvalue} (which is called the {\it principal eigenvalue}) of the differential operator $\cL_{a(\cd)}$ defined by the following:
\bel{L}\ba{ll}
\ns\ds\cL_{a(\cd)}\f(\cd)=-\nabla\cd\big(a(\cd)\nabla\f(\cd)\big),\qq\f(\cd)\in
\sD\big(\cL_{a(\cd)}\big),\\
\ns\ds\sD\big(\cL_{a(\cd)}\big)=\Big\{\f(\cd)\in W^{1,2}_0(\O)\bigm|\nabla\cd\big(a(\cd)\nabla\f(\cd)\big)\in L^2(\O)\Big\}.\ea\ee
In the above, $W^{1,2}_0(\O)$ is the usual Sobolev space which is the completion of $C^\infty_0(\O)$ (smooth functions with compact supports in $\O$) under the norm (see \cite{Adams-1975}):
$$\|\f(\cd)\|_{W^{1,2}_0(\O)}=\(\|\nabla\f(\cd)\|_2^2+\|\f(\cd)\|_2^2\)^{1\over2}<\infty.$$
It is known that with such a $\l_{a(\cd)}>0$, the following boundary value problem
$$\left\{\2n\ba{ll}
\ds-\nabla\cd\big(a(x)\nabla y_1(x)\big)=\l_{a(\cd)}y_1(x),\qq\hb{in }\O,\\
\ns\ds y_1\big|_{\pa\O}=0\ea\right.$$
admits a weak solution $y_1(\cd)\in W^{1,2}_0(\O)\setminus\{0\}$, which is called a {\it principal eigenfunction} of operator $\cL_{a(\cd)}$. Moreover,
\bel{l_A}\l_{a(\cd)}\equiv\inf_{y(\cd)\in
W^{1,2}_0(\O)\setminus\{0\}}{\ds\int_\O\lan a(x)\nabla y(x),\nabla
y(x)\ran dx\over\|y(\cd)\|^2_2}.\ee
From \cite{Gilbarg-Trudinger-1998}, Theorem 8.38, we know that the
multiplicity of $\l_{a(\cd)}$ is 1, and $y_1(\cd)\equiv
y_1(\cd\,;a(\cd))$ can be taken the unique eigenfunction such that
it is positive in $\O$ and normalized:
$$\|y_1(\cd)\|_2^2=\int_\O|y_1(x)|^2dx=1.$$
We call such a $y_1(\cd)$ the {\it normalized principal eigenfunction} of $\cL_{a(\cd)}$ and denote it by $y_{a(\cd)}$, indicating the dependence on $a(\cd)$. For convenience, we call $(\l_{a(\cd)},y_{a(\cd)})$ the {\it normalized principal eigen-pair} of the operator $\cL_{a(\cd)}$.
From (\ref{1.2}), we see that with such a $\l_{a(\cd)}$, one has
$$\|y(t,\cd\,;y_0(\cd))\|_2\les e^{-\l_{a(\cd)}t}\|y_0(\cd)\|_2,\qq t\ges0,\q\forall y_0(\cd)
\in L^2(\O)$$
and if the initial state $y_0(\cd)=y_{a(\cd)}$, then one has the equality:
$$y(t,x;y_{a(\cd)})=e^{-\l_{a(\cd)}t}y_{a(\cd)}(x),\qq(t,x)\in[0,\infty)\times\dbR^n,$$
which leads to
$$\|y(t,\cd\,;y_{a(\cd)})\|_2=e^{-\l_{a(\cd)}t}\|y_{a(\cd)}\|_2
=e^{-\l_{a(\cd)}t},\qq t\ges0.$$
Hence, one obtains the following representation:
$$\l_{a(\cd)}=\inf_{y_0(\cd)\in W^{1,2}_0(\O)\setminus\{0\}}\liminf_{t\to\infty}{1\over t}\log{\|y_0(\cd)\|_2\over\|y(t,\cd\,;y_0(\cd))\|_2}=-{1\over
t}\log\|y(t,\cd;y_{a(\cd)})\|_2,\qq\forall t>0.$$
Consequently, in some sense, $\l_{a(\cd)}$ is the smallest (uniform)
decay rate for the evolutionary map $y_0(\cd)\mapsto
y(t,\cd\,;y_0(\cd))$ (uniform with respect to the initial state $y_0(\cd)$).

\ms

On the other hand, since $a(\cd)$ is assumed to be uniformly
elliptic, the following Poincar\'e's inequality always holds:
\bel{P}\|y(\cd)\|_2^2\les C\int_\O\lan a(x)\nabla y(x),\nabla
y(x)\ran dx,\q\forall y(\cd)\in W^{1,2}_0(\O)\ee
for some constant $C>0$. Thus, the sharp constant $C$ that makes the
above true is given by the following
\bel{}\sup_{y(\cd)\in W_0^{1,2}(\O)\setminus\{0\}}{\|y(\cd)\|_2^2\over
\ds\int_\O\lan a(x)\nabla y(x),\nabla y(x)\ran
dx}\equiv{1\over\l_{a(\cd)}}<\infty.\ee
Therefore, the sharp Poincar\'e's inequality reads
\bel{P1}\|y(\cd)\|_2^2\les{1\over\l_{a(\cd)}}\int_\O\lan a(x)\nabla
y(x),\nabla y(x)\ran dx,\q\forall y(\cd)\in W^{1,2}_0(\O).\ee

If the diffusibility matrix $a(\cd)$ can be chosen from a given set $\sA$, which amounts to saying that the composite material/medium occupying $\O$ can be designed within a certain range, then we may try to minimize $\l_{a(\cd)}$ (preserving the temperature of the body in a certain fashion), or to maximize $\l_{a(\cd)}$ (cooling down the body as quick as possible, uniformly in the initial temperature distribution). In terms of Poincar\'e's inequality, the former means that we are looking for the sharp constant uniform for $a(\cd)\in\sA$, and the latter means that we are looking for the smallest possible sharp constant for some $a(\cd)\in\sA$.
\ms

Now, let $0<\m_0\les\m_1<\infty$ be given and let
\bel{M}M[\m_0,\m_1]=\Big\{A\in\dbS^n\bigm|\m_0I\les A\les\m_1I\Big\},\ee
where $\dbS^n$ is the set of all $(n\times n)$ symmetric matrices. Define
\bel{sM}\sM[\m_0,\m_1]=\Big\{a:\O\to M[\m_0,\m_1]\bigm|a(\cd)\hb{ is measurable}\Big\}
\subseteq L^\infty(\O;\dbS^n).\ee
Clearly, $M[\m_0,\m_1]$ is convex and compact in $\dbS^n$. Consequently, $\sM[\m_0,\m_1]$ is convex and closed in $L^p(\O;\dbS^n)$ (for any $p\in[1,\infty]$). From \eqref{l_A}, we see that
\bel{a<bar a}\l_{a(\cd)}\les\l_{\bar a(\cd)},\qq\forall a(\cd),\bar a(\cd)\in\sM[\m_0,\m_1],~a(\cd)\les\bar a(\cd),\ee
namely, the map $a(\cd)\mapsto\l_{a(\cd)}$ is monotone non-decreasing. In particular,
\bel{l<l<l}\m_0\l_I\les\l_{a(\cd)}\les\m_1\l_I,\qq\forall
a(\cd)\in\sM[\m_0,\m_1],\ee
where $\l_I$ is the principal eigenvalue of $-\D$ on $\O$, with the homogeneous Dirichlet
boundary condition. Consequently,
\bel{1.10}
\sup_{a(\cd)\in\sM[\m_0,\m_1]}\l_{a(\cd)}=\m_1\l_I,\qq\inf_{a(\cd)\in\sM[\m_0,\m_1]}\l_{a(\cd)}=\m_0\l_I.\ee
This shows that minimizing or maximizing $\l_{a(\cd)}$ over $\sM[\m_0,\m_1]$ is trivial. Now, we let $\sA$ be chosen with
\bel{1.11}\varnothing\ne\sA\subsetneq\sM[\m_0,\m_1],\qq0<\m_0\les\m_1<\infty,\ee
and pose the following two problems.

\ms

\bf Problem ($\bar\L(\sA)$). \rm Find an $\bar a(\cd)\in\sA$ such that
\bel{max}\l_{\bar a(\cd)}=\sup_{a(\cd)\in\sA}\l_{a(\cd)}.\ee

\ms

\bf Problem ($\underline\L(\sA)$). \rm Find an $\underline a(\cd)\in\sA$ such that
\bel{min}\l_{\underline a(\cd)}=\inf_{a(\cd)\in\sA}\l_{a(\cd)}.\ee

\ms

Some general results will be presented concerning the above two problems in the next section.

\ms

Further, to obtain finer results, we will concentrate on a more specific case which we now describe. Fix two different matrices $A_0,A_1\in M[\m_0,\m_1]$ with $0<\m_0\les\m_1<\infty$. We define
\bel{A(r)}A(r)=(1-r)A_0+rA_1,\qq r\in[0,1],\ee
and for some $0\les\a\les\b\les1$, let
\bel{cU}\sU[\a,\b]=\Big\{\chi_{_{\O_1}}(\cd)\bigm|\O_1\subseteq\O\hb{ is measurable, }
\a|\O|\les|\O_1|\les\b|\O|\Big\},\ee
where $\chi_{_{\O_1}}(\cd)$ is the characteristic function of $\O_1$ and $|\O_1|$ is the Lebesgue measure of measurable set $\O_1$. Then for any $u(\cd)\equiv\chi_{_{\O_1}}(\cd)\in\sU[\a,\b]$, we have
$$A(u(\cd))=[1-\chi_{_{\O_1}}(\cd)]A_0+\chi_{_{\O_1}}(\cd)A_1=A_0+\chi_{_{\O_1}}(\cd)(A_1-A_0).$$
Denote
\bel{sA[a,b]}\sA[\a,\b]=A\big(\sU[\a,\b]\big)\1n=\1n\Big\{A_0+\chi_{_{\O_1}}(\cd)(A_1-A_0)
\bigm|\O_1\hb{ measurable, }~\a|\O|\1n\les\1n|\O_1|\1n\les\1n\b|\O|\Big\}.\ee
Note that the set $\sA[\a,\b]$ is non-convex (unless $A_0=A_1$ which is excluded). For the heat conduct problem, with $a(\cd)=A(u(\cd))$ for some $u(\cd)=\chi_{_{\O_1}}(\cd)\in\sU[\a,\b]$, it means that two media occupy the domain, the one with conductivity matrix $A_1$ occupies
$\O_1$ and the other with conductivity matrix $A_0$ occupies $\O\setminus\O_1$. The corresponding principal eigenvalue and the corresponding (unique) normalized principal eigenfunction are denoted by
$$\l_{u(\cd)}\equiv\l_{A_0+u(\cd)(A_1-A_0)},\qq y_{u(\cd)}\equiv y_{A_0+u(\cd)(A_1-A_0)},$$
and the following holds:
\bel{}\left\{\2n\ba{ll}
\ds-\nabla\cd\(\big[A_0+u(x)(A_1-A_0)\big]\nabla y_{u(\cd)}(x)\)=\l_{u(\cd)}y_{u(\cd)}(x),\qq x\in\O,\\
\ns\ds y_{u(\cd)}\big|_{\pa\O}=0.\ea\right.\ee
Then we can pose the following problems.

\ms

\bf Problem ($\bar\L[\a,\b]$). \rm For given $0\les\a\les\b\les1$, find a $\bar u(\cd)\in\sU[\a,\b]$ such that
\bel{sup}\l_{\bar u(\cd)}=\sup_{u(\cd)\in\sU[\a,\b]}\l_{u(\cd)}.\ee

\ms

\bf Problem ($\underline\L[\a,\b]$). \rm For given $0\les\a\les\b\les1$, find a $\underline u(\cd)\in\sU[\a,\b]$ such that
\bel{inf}\l_{\underline u(\cd)}=\inf_{u(\cd)\in\sU[\a,\b]}\l_{u(\cd)}.\ee

\ms

Any $\bar u(\cd)\in\sU[\a,\b]$ (resp. $\underline u(\cd)\in\sU[\a,\b]$) satisfying (\ref{sup}) (resp.$\,$(\ref{inf})) is called an {\it optimal control} of Problem ($\bar\L[\a,\b]$)
(resp. Problem ($\underline\L[\a,\b]$)).

\ms

Note that if $0=\a<\b\les1$, and, say, $A_0\les A_1$, due to a fact similar to (\ref{l<l<l})--(\ref{1.10}), one has
\bel{inf**}\inf_{u(\cd)\in\sU[0,\b]}\l(u(\cd))=\l_{A_0},\ee
making Problem ($\underline\L[0,\b]$) trivial; and likewise, if $0<\a\les\b=1$, and still let $A_0\les A_1$, then
\bel{sup*}\sup_{u(\cd)\in\sU[\a,1]}\l(u(\cd))=\l_{A_1},\ee
making Problem ($\bar\L[\a,1]$) trivial. To avoid such situations,
in what follows, we will assume the following:
\bel{ab}\left\{\1n\ba{ll}
\ds\hb{either }0<\a\les\b<1,\q\hb{no additional restrictions on $A_0,A_1$},\\
\ns\ds\hb{or $\a=0$, $\b=1$,~~and~~~neither $A_0\les A_1$, nor $A_1\les A_0$ holds}.\ea\right.\ee
Note that when $0<\a\les\b<1$, even if, say, $A_1=2A_0>A_0$, the location/shape of the optimal $\O_1$ (if it exists, which is non-empty and not equal to $\O$) is not obvious. On the other hand, in the case that $A_0$ and $A_1$ are not comparable, one expects that neither $u_0(\cd)=0$ nor $u_1(\cd)=1$ is optimal.

\ms

For either case in \eqref{ab}, $\sU[\a,\b]$ is not convex. Hence, the existence of optimal controls for Problems ($\bar\L[\a,\b]$) and ($\underline\L[\a,\b]$) is not guaranteed, in general. To study these problems, we will introduce suitable relaxed problems for which the relaxed optimal controls will exist. Some necessary conditions for relaxed optimal controls will then be established, and illustrative examples will be presented as well.

\ms

Some studies on optimization of the principal eigenvalue for elliptic operators can be found in the book by Henrot \cite{Henrot-2006} (see also the references cited therein). The case studied in \cite{Henrot-2006} was isotropic, namely, the diffusion matrix $a(x)=\si(x)I$, for some scalar function $\si(\cd)$. Moreover, even for that case, only a maximization problem was considered. For other relevant works, here is a partial list of references: \cite{Cox-McLaughlin-1990, Cox-Kawohl-Uhlig-1999, Belhachmi-Bucur-Buttazzo-Sac-Epee-2006, Caffarelli-Lin-2007, Munch-Pefregel-Periago-2008, Cuccu-Emamizadeh-Porru-2008, Casado-Diaz-2015, Bucur-Buttazzo-Nitsch-2017}.

\ms

The rest of this paper is organized as follows. Section 2 will be devoted to some general considerations of the problems that we are interested in. In Section 3, a convexification of maximization problem is investigated. In Section 4, a relaxation of minimization problem in terms of the so-called $H$-convergence will be studied. A detailed example is worked out in Section 5.
Finally some remarks are collected in Section 6.

\section{The Principal Eigen Map and Its Properties}

We fix a bounded Lipschitz domain $\O\subset\dbR^n$, i.e., $\O\subseteq\dbR^n$ is a bounded domain with a Lipschitz boundary $\pa\O$, and constants $0<\m_0\les\m_1<\infty$. For any $a(\cd)\in\sM[\m_0,\m_1]$, recall that $\cL_{a(\cd)}$ is an elliptic operator defined by (\ref{L}), and $(\l_{a(\cd)},y_{a(\cd)})$ is the normalized principal eigen-pair of $\cL_{a(\cd)}$. Thus,
\bel{2.1}\left\{\2n\ba{ll}
\ds-\nabla\cd\big(a(x)\nabla
y_{a(\cd)}(x)\big)=\l_{a(\cd)}y_{a(\cd)}(x),\qq\hb{in }\O,\\
\ns\ds y_{a(\cd)}\big|_{\pa\O}=0,\ea\right.\ee
with $y_{a(\cd)}(x)>0$ for $x\in\O$, and $\|y_{a(\cd)}(\cd)\|_2=1$. Further, it is known that
\bel{2.2}\l_{a(\cd)}=\min_{y(\cd)\in W^{1,2}_0(\O),\|y(\cd)\|_2=1}\int_\O\lan a(x)\nabla y(x),\nabla y(x)\ran dx=\int_\O\lan a(x)\nabla y_{a(\cd)}(x),\nabla
y_{a(\cd)}(x)\ran dx.\ee
Define $\L:\sM[\m_0,\m_1]\to\dbR\times W^{1,2}_0(\O)$ by
\bel{F}\L(a(\cd))=(\l_{a(\cd)},y_{a(\cd)}),\qq\forall a(\cd)\in\sM[\m_0,\m_1]\ee
and call it the {\it principal eigen map} of the operator $\cL_{a(\cd)}$. We first recall the following result found in \cite{Huska-Polacik-Safonov-2005}.

\bp{Proposition 2.4} \sl Let $a(\cd)\in\sM[\m_0,\m_1]$. Then there exists a constant $\g=\g(\m_0,\m_1;\O)>0$ such that
\bel{2.9}\l-\l_{a(\cd)}\ges\g,\qq\ee
where $\l$ is any eigenvalue of $\cL_{a(\cd)}$ different from $\l_{a(\cd)}$. Consequently, one has:
\bel{l-l_a>g}\ba{ll}
\ns\ds\int_\O\2n\lan a(x)\nabla y(x),\nabla y(x)\ran dx\1n\ges\1n\big(\l_{a(\cd)}\1n+\1n\g\big)\,\|y(\cd)\|_2^2,\q~\forall y(\cd)\1n\in\1n W^{1,2}_0(\O),\;\int_\O y(x)y_{a(\cd)}(x)dx=0.\ea\ee

\ep

Next, we present the following simple result which will be useful below.

\bp{Proposition 2.1} \sl The map $a(\cd)\mapsto\l_{a(\cd)}$ is concave on $\sM[\m_0,\m_1]$. Consequently, this map is Lipschitz on $L^\infty(\O;\dbS^n)$, i.e., for some $L>0$,
\bel{|l-l|<|a-a|_infty}|\l_{a(\cd)}-\l_{\bar a(\cd)}|\les L\|a(\cd)-\bar a(\cd)\|_\infty,\qq\forall a(\cd),\bar a(\cd)\in L^\infty(\O;\dbS^n).\ee

\ep

\it Proof. \rm For any $a_0(\cd),a_1(\cd)\in\sM[\m_0,\m_1]$ and any $\g\in(0,1)$,
$$\ba{ll}
\ns\ds\l_{(1-\g)a_0(\cd)+\rho a_1(\cd)}=\min_{y(\cd)\in W_0^{1,2}(\O)\atop \|y(\cd)\|_2=1}\int_\O\ip{\big[(1-\g)a_0(x)+\g
a_1(x)\big]\nabla y(x),\nabla y(x)}\, dx\\
\ns\ds\ges(1-\g)\2n\min_{{y(\cd)\in W_0^{1,2}(\O)}\atop{\|y(\cd)\|_2=1}}\int_\O\lan a_0(x)\nabla
y(x),\nabla y(x)\ran dx\1n+\1n\g\2n\min_{{y(\cd)\in W^{1,2}_0(\O)}\atop{\|y(\cd)\|_2=1}}\int_\O\lan a_1(x)\nabla
y(x),\nabla y(x)\ran dx\\
\ns\ds\ges(1-\g)\l_{a_0(\cd)}+\g\l_{a_1(\cd)}.\ea$$
This proves the concavity of the map $a(\cd)\mapsto\l_{a(\cd)}$. Then the Lipschitz continuity follows from a standard argument (see, for example, \cite{Li-Yong-1995}, p.235, Lemma 2.8). \endpf

\ms

We will see that the map $a(\cd)\mapsto\l_{a(\cd)}$ is not strictly convex.

\ms

The following theorem is due to Gallouet--Monier \cite{Gallouet-Monier-2000}, which is an extension of a result by Meyers \cite{Meyers-1963}.

\bt{Meyers-Gallouet-Monier} \sl Let $\O\subseteq\dbR^n$ be a bounded Lipschitz domain and $0<\m_0\les\m_1$. Then there exists a $p_0>2$, only depending on $\O$, such that for any $2<p<p_0$, $a(\cd)\in\sM[\m_0,\m_1]$,  $\Bf(\cd)\in L^p(\O)$ and $h(\cd)\in L^{np\over n+p}(\O)$, the following problem:
$$\left\{\2n\ba{ll}
\ns\ds-\nabla\cd\(a(x)\nabla y(x)\)=\nabla\cd\Bf(x)+h(x),\qq
x\in\O,\\
\ns\ds y\big|_{\pa\O}=0\ea\right.$$
admits a unique weak solution $y(\cd)\in W^{1,p}_0(\O)$, and the
following estimate holds:
\bel{|y'|_p}\|\nabla y(\cd)\|_p\les C\(\|\Bf(\cd)\|_p+\|h(\cd)\|_{np\over n+p}\).\ee
Hereafter, $C>0$ represents a generic constant which could be different from line to line.

\et

Note that in the above result, $p>2$, which will play a crucial role below. The following result is concerned with the principal eigen map $\L$.

\bt{Theorem 2.5} \sl Let $0<\m_0\les\m_1<\infty$ and $\O\subseteq\dbR^n$ be a Lipschitz domain. Then there exists a $p>2$ such that the principal eigen map $\L$ is Lipschitz continuous in the following sense:
\bel{l-l}|\l_{a_1(\cd)}-\l_{a_2(\cd)}|\les C\|a_1(\cd)-a_2(\cd)\|_{p\over p-2},\qq\forall a_1(\cd),a_2(\cd)\in\sM[\m_0,\m_1],\ee
and for any $\bar p\in[2,p)$,
\bel{|yA-yA|}\|y_{a_1(\cd)}-y_{a_2(\cd)}\|_{W^{1,\bar p}(\O)}\les
C\|a_1(\cd)-a_2(\cd)\|_{{p\bar p\over p-\bar p}},\qq\forall a_1(\cd),a_2(\cd)\in\sM[\m_0,\m_1].\ee

\et

Note that \eqref{l-l} is an improvement of \eqref{|l-l|<|a-a|_infty}, thanks to the existence of a $p>2$ so that \eqref{|y'|_p} holds.

\ms

\it Proof. \rm First of all, for any $a(\cd)\in\sM[\m_0,\m_1]$, we recall that
$$\m_0\|\nabla y_{a(\cd)}\|_2^2\les\int_\O\lan a(x)\nabla
y_{a(\cd)}(x),\nabla y_{a(\cd)}(x)\ran dx=\l_{a(\cd)}\les\m_1\l_I$$
and by Sobolev embedding theorem (\cite{Adams-1975}), for any
$n\ges2$, we have
$$W^{1,2}(\O)\hookrightarrow L^r(\O),\qq r<{2n\over n-2}$$
with the convention that ${2n\over n-2}=\infty$ when $n=2$.

\ms

Now, let $p_0>0$ be the number in  Theorem \ref{Meyers-Gallouet-Monier} and chose $2<p<\min\big(p_0,{2n\over  n-2}\big)$. We have $p(n-2)<2n$, which leads to ${np\over n+p}<2$. Hence, by Theorem \ref{Meyers-Gallouet-Monier}, regarding $\l_{a(\cd)}y_{a(\cd)}$ as a nonhomogeneous term on the right-hand side of the equation (\ref{2.1}),
we have (noting $\|y_{a(\cd)}\|_2=1$)
\bel{|y_a|}\|y_{a(\cd)}\|_{W^{1,p}(\O)}\les
C_p\l_{a(\cd)}\|y_{a(\cd)}\|_{np\over n+p}\les C_p\l_{a(\cd)}|\O|^{{1\over n}+{1\over p}-{1\over 2}}=\bar C_p.\ee
Here $\bar C_p$ is an absolute constant, uniform in $a(\cd)\in\sM[\m_0,\m_1]$. Let
$$\sY^p(\m_0,\m_1)=\big\{y(\cd)\in W^{1,p}_0(\O)\bigm|\|y(\cd)\|_2=1,~
\|y(\cd)\|_{W^{1,p}(\O)}\les\bar C_p,~y(x)\ges0,~\ae\big\}.$$
Then $y_{a(\cd)}\in\sY^p(\m_0,\m_1)\subseteq W^{1,2}_0(\O)$. Consequently, for any $a_1(\cd),a_2(\cd)\in\sM[\m_0,\m_1]$, denote
$$\l_i=\l_{a_i(\cd)},\q y_i(\cd)=y_{a_i(\cd)}(\cd),\qq i=1,2,$$
and assume that $\l_1\ges\l_2$. We have
$$\ba{ll}
\ns\ds\l_1-\l_2=\l_1-\int_\O\lan a_2(x)\nabla y_2(x),\nabla y_2(x)\ran\, dx\\
\ns\ds\les\int_\O\lan a_1(x)\nabla y_2(x),\nabla y_2(x)\ran
dx-\int_\O\lan a_2(x)\nabla y_2(x),\nabla y_2(x)\ran\, dx\\
%
\ns\ds\les\|a_1(\cd)-a_2(\cd)\|_{p\over p-2}\|\nabla
y_2(\cd)\|_p^2\les\bar C_p^2\|a_1(\cd)-a_2(\cd)\|_{p\over p-2},\ea$$
proving \eqref{l-l}. Next, let
$$\ba{ll}
\ns\ds \a=\int_\O y_1(\xi)y_2(\xi)d\xi,\q \\
\ns\ds y_{12}(x)=y_1(x)-\a y_2(x) \q y_{21}(x)=y_2(x)-\a y_1(x).\ea$$
Then $\a\ges0$, and
$$\ba{ll}
\ns\ds\int_\O y_{12}(x)y_2(x)dx=\int_\O y_{21}(x)y_1(x)dx=0,\qq\|y_{12}(\cd)\|_2^2=\|y_{21}(\cd)\|_2^2=1-\a^2,\\ [2mm]
\ns\ds(1+\a)\big(y_2(x)-y_1(x)\big)=y_{21}(x)-y_{12}(x).\ea$$
Hence,
\bel{N2001}\|y_2(\cd)-y_1(\cd)\|_2^2\les  \|y_{21}(\cd)-y_{12}(\cd)\|^2=2\|y_{21}(\cd)\|_2^2+2\|y_{12}(\cd)\|^2=4(1-\a^2).\ee
By \eqref{l-l_a>g}, one has
$$\ba{ll}
\ns\ds(1-\a^2)\g=\l_1 \a^2+(\l_1+\g)(1-\a^2)-\l_1\\
\ns\ds\les\a^2\int_\O\lan a_1(x)\nabla y_1(x),\nabla y_1(x)\ran dx+
 \int_\O\lan a_1(x)\nabla y_{21}(x),\nabla y_{21}(x)\ran dx-\l_1\\
\ns\ds=\int_\O\lan a_1(x)\nabla y_2(x),\nabla y_2(x)\ran dx-\l_1\\
\ns\ds=\l_2-\l_1+\int_\O\lan(a_1(x)-a_2(x))\nabla y_1(x),\nabla y_2(x)\ran dx\\
\ns\ds\qq+\int_\O\lan(a_1(x)-a_2(x)\big)\big(\na y_2(x)-\nabla y_1(x)\big),\nabla y_2(x)\ran dx\\
\ns\ds\les(1-\a)(\l_2-\l_1) +\|a_1(\cd)-a_2(\cd)\|_{2p\over p-2}\|\nabla y_1(\cd)-\nabla y_2(\cd)\|_2\|\nabla y_2(\cd)\|_p.\ea$$
Consequently,
$$(1-\a^2)\g\les C\|a_1(\cd)-a_2(\cd)\|_{2p\over p-2}\|\nabla y_1(\cd)-\nabla y_2(\cd)\|_2.$$
Combining the above with \eqref{N2001}, we have
\bel{y_1-y_2}\|y_1(\cd)-y_2(\cd)\|_2^2\les C\|a_1(\cd)-a_2(\cd)\|_{2p\over p-2}\|\nabla y_1(\cd)-\nabla y_2(\cd)\|_2.\ee
On the other hand,
$$\left\{\2n\ba{ll}
\ds-\nabla\cd\big(a_1(x)\nabla[y_1(x)-y_2(x)]\big)=\l_1y_1-\l_2y_2
+\nabla\cd\big([a_1(x)-a_2(x)]\nabla y_2(x)\big),\\
\ns\ds\big(y_1-y_2\big)\big|_{\pa\O}=0.\ea\right.$$
It follows from Theorem \ref{Meyers-Gallouet-Monier} that (with $2\les\bar p<p$)
\bel{N2003}\ba{ll}
\ns\ds\|y_1(\cd)-y_2(\cd)\|_{W^{1,\bar p}(\O)}\les C\|(a_1(\cd)-a_2(\cd))\nabla y_2\|_{\bar
p}+C\|\l_1 y_1(\cd)-\l_2 y_2(\cd)\|_{n\bar p\over n+\bar p}\\
\ns\ds\les C\|a_1(\cd)-a_2(\cd)\|_{p\bar p\over p-\bar p}\|\nabla y_2(\cd)\|_p+C|\O|^{{1\over n}+{1\over\bar p}-{1\over 2}}\|\l_1 y_1(\cd)-\l_2 y_2(\cd)\|_2\\
\ns\ds\les C\bar C_p\|a_1(\cd)-a_2(\cd)\|_{p\bar p\over p-\bar p}+C|\O|^{{1\over n}+{1\over\bar p}-{1\over 2}}\,\Big(|\l_1-\l_2|+ \m_1 \|y_1(\cd)-y_2(\cd)\|_2\Big)\\
\ns\ds\les C\(\|a_1(\cd)-a_2(\cd)\|_{p\bar p\over p-\bar p}+\|y_1(\cd)-y_2(\cd)\|_2\).\ea\ee
Applying \eqref{y_1-y_2}, we have
$$\ba{ll}
\ns\ds\|y_1(\cd)-y_2(\cd)\|_{W^{1,\bar p}(\O)}\les C\(\|a_1(\cd)-a_2(\cd)\|_{p\bar p\over p-\bar p}+\|a_1(\cd)-a_2(\cd)\|_{2p\over p-2}^{1\over2}\|y_1(\cd)-y_2(\cd)\|_{W^{1,2}(\O)}^{1\over2}\)\\
\ns\ds\les C\(\|a_1(\cd)-a_2(\cd)\|_{p\bar p\over p-\bar p}+\|a_1(\cd)-a_2(\cd)\|_{p\bar p\over p-\bar p}^{1\over2}\|y_1(\cd)-y_2(\cd)\|_{W^{1,\bar p}(\O)}^{1\over2}\)\\
\ns\ds\les C\|a_1(\cd)-a_2(\cd)\|_{p\bar p\over p-\bar p}+{1\over2}\|y_1(\cd)-y_2(\cd)\|_{W^{1,\bar p}(\O)}.\ea$$
Then \eqref{|yA-yA|} follows. \endpf



We now look at the directional differentiability of the principal eigen map
$\L$.

\bp{Proposition 2.8} \sl Let $\O$ be a bounded Lipschitz domain and $a(\cd),\bar a(\cd)\in\sM[\m_0,\m_1]$. Then the directional derivative of the eigen map $\L$ at $\bar a(\cd)$ in the direction of $a(\cd)-\bar a(\cd)$ is given by
\bel{E215}\L'\big(\bar a(\cd);a(\cd)-\bar a(\cd)\big)\equiv\lim_{\e\da0}{\L(\bar a(\cd)+\e[a(\cd)-\bar a(\cd)])-\L(\bar a(\cd))\over\e}=\big(\l',y'(\cd)\big),\ee
where $\l'$ is given by the following:
\bel{hl}\l'=\int_\O\lan[a(x)-\bar a(x)]\nabla y_{\bar
a(\cd)}(x),\nabla y_{\bar a(\cd)}(x)\ran dx,\ee
and $y'(\cd)$ is the weak solution to the following:
\bel{limit equation}\left\{\2n\ba{ll}
\ds-\nabla\cd\big(\bar a(x)\nabla
y'(x)\big)=\l_{\bar a(\cd)}y'(x)+\l'\,y_{\bar a(\cd)}(x)+\nabla\cd\big([a(x)-\bar
a(x)]\nabla y_{\bar a(\cd)}(x)\big),\\
\ns\ds y'\big|_{\pa\O}=0.\ea\right.\ee

\ep

\it Proof. \rm Let $\e\in (0,1)$. Denote
$$a_\e(\cd)=\bar a(\cd)+\e[a(\cd)-\bar a(\cd)]\in\sM(\m_0,\m_1),$$
and
$$\L(\bar a(\cd))=\big(\l_{\bar a(\cd)},y_{\bar a(\cd)}\big)\equiv(\bar\l,\bar
y(\cd)),\qq\L(a_\e(\cd))=(\l_{a_\e(\cd)},y_{a_\e(\cd)})\equiv(\l_\e,y_\e(\cd)).$$
Let
$$\l'_\e={\l_\e-\bar\l\over\e},\qq y'_\e(\cd)={y_\e(\cd)-\bar y(\cd)\over\e}.$$
Note that
$$\ba{ll}
\ns\ds{a_\e(x)\nabla y_\e(x)-\bar a(x)\nabla\bar y(x)\over\e}={\big[\bar a(x)+\e\big(a(x)-\bar a(x)\big)\big]\nabla y_\e(x)-\bar a(x)\nabla\bar y(x)\over\e}\\
\ns\ds\qq\qq\qq\qq\qq\qq=\bar a(x)\nabla y'_\e(x)+\big[a(x)-\bar a(x)\big]\nabla y_\e(x),\ea$$
and
$${\l_\e y_\e(x)-\bar\l\bar y(x)\over\e}={\l_\e[y_\e(x)-\bar y(x)]+(\l_\e-\bar\l)\bar y(x)\over\e}
=\l_\e y'_\e(x)+\l'_\e\bar y(x).$$
Hence, $(\l_\e',y_\e'(\cd))$ satisfies
\bel{2.21}\left\{\2n\ba{ll}
\ds-\nabla\cd\big(\bar a(x)\nabla y'_\e(x)\big)=\l_\e y'_\e(x)+\l'_\e\bar y(x)+\nabla\cd\big([a(x)-\bar a(x)]\nabla y_\e(x)\big),\\
\ns\ds y'_\e\big|_{\pa\O}=0.\ea\right.\ee
Since $a(\cd)\mapsto\l_{a(\cd)}$ is Lipschitz, we have
$$|\l'_\e|={|\l_\e-\bar\l|\over\e}\les L\|a(\cd)-\bar a(\cd)\|_\infty.$$

From the Lipschitz continuity of the principal eigen map $\L$ (see Theorem \ref{Theorem 2.5}), we have
\bel{|y'|<C}|\l'_\e|+\|y'_\e(\cd)\|_{W^{1,\bar p}(\O)}\les C,\qq\forall\e>0,\ee
with $2\les\bar p<p$. Then we may assume that
$$\l'_\e\to\l',\q\hb{and}\q y'_\e(\cd)\to y'(\cd),\q\hb{weakly in $W^{1,\bar p}_0(\O)$}.$$
Also, \eqref{|y'|<C} leads to
$$\l_\e\to\bar\l,\q\hb{and}\q y_\e(\cd)\to\bar y(\cd),\q\hb{strongly in $W^{1,\bar p}(\O)$}.$$
Consequently, passing to the limit in \eqref{2.21}, one gets
\bel{limit equation1}\left\{\2n\ba{ll}
\ds-\nabla\cd\big(\bar a(x)\nabla
y'(x)\big)=\bar\l y'(x)+\l'\bar y(x)+\nabla\cd\big([a(x)-\bar
a(x)]\nabla\bar y(x)\big),\\
\ns\ds y'\big|_{\pa\O}=0,\ea\right.\ee
and
$$\ba{ll}
\ns\ds\l'=\l'\|\bar y(\cd)\|_2^2=\int_\O\bar y(x)\(-\nabla\cd\big(\bar a(x)\nabla y'(x)\big)-\bar\l y'(x)-\nabla\cd\big([a(x)-\bar a(x)]\nabla\bar y(x)\big)\)dx\\
\ns\ds\q=\int_\O\(\lan\bar a(x)\nabla\bar y(x),\nabla y'(x)\ran-\bar\l\bar y(x)y'(x)+\lan\nabla\bar y(x),[a(x)-\bar a(x)]\nabla\bar y(x)\ran\)dx\\
\ns\ds\q=\int_\O\lan\,[a(x)-\bar a(x)]\nabla\bar y(x),\nabla\bar y(x)\ran dx.\ea$$
This proves our conclusions. \endpf

The following gives some direct consequences of the above general results.

\bc{Corollary 2.7} \sl Let $\O$ be a bounded Lipschitz domain. Let $0<\m_0\les\m_1<\infty$, and $\sA\subseteq\sM[\m_0,\m_1]$ be convex and closed in $L^1(\O;\dbS^n)$. Then the following conclusions are true:

\ms

{\rm(i)} Problem ($\bar\L(\sA)$) admits an optimal solution. Further, $\bar a(\cd)\in
\sM[\m_0,\m_1]$ is an optimal solution to Problem ($\bar\L(\sA)$) if and only if
\bel{A-A<0}\int_\O\lan[a(x)-\bar a(x)]\nabla y_{\bar a(\cd)}(x),\nabla y_{\bar a(\cd)}(x)\ran dx\les0,\qq\forall a(\cd)\in\sA.\ee

{\rm(ii)} If $\underline a(\cd)\in\sM[\m_0,\m_1]$ is an optimal solution
to Problem ($\underline\L(\sA)$). Then
\bel{A-A>0}\int_\O\lan[a(x)-\underline a(x)]\nabla y_{\underline
a(\cd)}(x),\nabla y_{\underline a(\cd)}(x)\ran dx\ges0,\qq\forall a(\cd)\in\sA.\ee

\ec

\it Proof. \rm (i) From Proposition \ref{Proposition 2.1} and Theorem \ref{Theorem 2.5}, we know that $a(\cd)\mapsto\l_{a(\cd)}$ is concave and continuous from $\sM[\m_0,\m_1]$ to $\dbR$. Hence, a standard argument involving Mazur's Theorem applies to get the existence of an optimal solution
to Problem ($\bar\L(\sA)$).

\ms

Next, if $\bar a(\cd)\in\sA$ is a maximum of $a(\cd)\mapsto\l_{a(\cd)}$, then for any
$a(\cd)\in\sA$, making use of the convexity of $\sA$, we have
$$0\ges\lim_{\e\da0}{\l_{\bar a(\cd)+\e[a(\cd)-\bar a(\cd)]}-\l_{\bar a(\cd)}\over\e}
=\l'=\int_\O\lan[a(x)-\bar a(x)]\nabla y_{\bar a(\cd)}(x),\nabla y_{\bar a(\cd)}(x)\ran dx.$$
This gives \eqref{A-A<0}. Conversely, suppose \eqref{A-A<0} holds. Then by the concavity of $a(\cd)\mapsto\l_{a(\cd)}$, we see that for any $a(\cd)\in\sA$, $\e\mapsto\l_{\bar a(\cd)+\e
[a(\cd)-\bar a(\cd)]}$ is concave as well. Thus, $\e\mapsto{d\over d\e}\l_{\bar a(\cd)+\e[a(\cd)-
\bar a(\cd)]}$ is non-increasing. Consequently,
$$\ba{ll}
\ns\ds\l_{a(\cd)}-\l_{\bar a(\cd)}={1\over\e}\[(1-\e)\l_{\bar a(\cd)}+\e\l_{a(\cd)}
-\l_{\bar a(\cd)}\]\les{1\over\e}\[\l_{\bar a(\cd)+\e[a(\cd)-\bar a(\cd)]}
-\l_{\bar a(\cd)}\]\\
\ns\ds\les\[{d\over d\e}\l_{\bar a(\cd)+\e[a(\cd)-\bar a(\cd)]}\]\Big|_{\e=0}
=\int_\O\lan[a(x)-\bar a(x)]\nabla y_{\bar a(\cd)}(x),\nabla y_{\bar a(\cd)}(x)\ran dx\les0.\ea$$
Hence, $\bar a(\cd)$ is a maximum of $\l_{a(\cd)}$ over $\sA$.

\ms

(ii) From Proposition \ref{Proposition 2.8}, we see that
$$\l_{\underline a(\cd)+\e[a(\cd)-\underline a(\cd)]}=\l_{\underline
a(\cd)}+\e\int_\O\lan[a(x)-\underline a(x)]\nabla y_{\underline
a(\cd)}(x),\nabla y_{\underline a(\cd)}(x)\ran dx+o(\e).$$
Hence, if $\underline a(\cd)\in\sA$ is a solution to Problem ($\underline\L(\sA)$), then \eqref{A-A>0} holds. \endpf

\ms

We note that Corollary \ref{Corollary 2.7} part (i) gives the existence and characterization of optimal solutions to Problem ($\bar\L(\sA)$), thanks to the concavity of the map $a(\cd)\mapsto\l_{a(\cd)}$. Whereas, part (ii) of
Corollary \ref{Corollary 2.7} only gives a necessary condition for a possible
solution of Problem ($\underline\L(\sA)$), and no existence of optimal solution is guaranteed.

\section{A Convexification of Problem $(\bar\L[\a,\b])$}

Let us return to Problem $(\bar\L[\a,\b])$. Since $\sU[\a,\b]$ is not convex, the existence of optimal solution is not guaranteed. In this section, we consider a convexification of Problem $(\bar\L[\a,\b])$.

\ms

For $0\les\a\les\b\les1$, we introduce the following:
\bel{Si[a,b]}\Si[\a,\b]=\Big\{\si:\O\to[0,1]\bigm|\si(\cd)\hb{ measurable, }\a|\O|\les\int_\O\si(x)dx\les\b|\O|\Big\},\ee
which is convex and closed in $L^1(\O)$. Recalling $\sU[\a,\b]$ defined by \eqref{cU}, one has
\bel{coU}\cl{\co\{\sU[\a,\b]\}}^{L^1(\O;\dbR)}=\Si[\a,\b],\ee
where the left hand side of the above is the closed convex hull of $\sU[\a,\b]$ in $L^1(\O;\dbR)$. Now, for given $0\les\a\les\b\les1$, and $A_0,A_1\in M[\m_0,\m_1]$ with $0<\m_0\les\m_1<\infty$ such that \eqref{ab} holds, with $A(\cd)$ defined by \eqref{A(r)}, one sees that
$$A\big(\Si[\a,\b]\big)=\Big\{A_0+\si(\cd)(A_1-A_0)\bigm|\si(\cd)\in\Si[\a,\b]\Big\}$$
is convex and closed in $L^1(\O;\dbR)$. For any $\si(\cd)\in\Si[\a,\b]$, we consider the following state equation
\bel{N3002}\left\{\2n\ba{ll}
\ds-\nabla\cd\big(A(\si(x))\nabla y(x)\big)=\l
y(x),\qq x\in\O,\\
\ns\ds y\big|_{\pa\O}=0.\ea\right.\ee
Denote the corresponding principal eigenvalue and normalized principal
eigenfunction by
$$\l_{\si(\cd)}\equiv\l_{A_0+\si(\cd)(A_1-A_0)},\qq y_{\si(\cd)}=y_{A_0+\si(\cd)(A_1-A_0)}(\cd),$$
respectively. We pose the following convexified problem.

\ms

\bf Problem ($\bar\L^c[\a,\b]$). \rm Let $0\les\a\les\b\les1$, $0<\m_0\les\m_1<\infty$, and $A_0,A_1\in M[\m_0,\m_1]$ satisfying (\ref{ab}). Find a $\bar\si(\cd)\in\Si[\a,\b]$ such that
\bel{l-bar-si}\l_{\bar\si(\cd)}=\sup_{\si(\cd)\in\Si[\a,\b]}\l_{\si(\cd)}.\ee

\ms

Any $\bar\si(\cd)\in\Si[\a,\b]$ satisfying \eqref{l-bar-si} is called an {\it optimal control} of Problem ($\bar\L^c[\a,\b]$), which is also called an {\it optimal convexified relaxed control} of Problem ($\bar\L[\a,\b]$). The superscript ``$c$'' in $\bar\L^c[\a,\b]$ indicates the ``convexification''. Note that if $\bar\si(\cd)$ is an optimal control of Problem ($\bar\L^c[\a,\b]$) and
\bel{0<si<1}\big|\big(0<\bar\si(\cd)<1\big)\big|=0,\ee
where
$$\big(0<\bar\si(\cd)<1\big)=\big\{x\in\O\bigm|0<\bar\si(x)<1\big\},$$
then
$$\bar u(\cd)=\chi_{\O_1}(\cd),\qq\O_1=(\bar\si(\cd)=1)$$
is an optimal control of Problem ($\bar\L[\a,\b]$). The following result gives the necessary conditions when \eqref{0<si<1} fails.

\bt{Theorem 3.1} \sl Problem {\rm($\bar\L^c[\a,\b]$)} admits an optimal control $\bar\si(\cd)\in\Si[\a,\b]$. Suppose
\bel{(0<si<1)>0}|(0<\bar\si(\cd)<1)|>0,\ee
and $\bar y(\cd)$ is the corresponding optimal state. Then
\bel{3.3A}\left\{\2n\ba{ll}
\ds\lan(A_1-A_0)\nabla\bar y(x),\nabla\bar y(x)\ran\equiv C,\qq\ae
x\in\big(0<\bar\si(\cd)<1\big),\\
\ns\ds\lan(A_1-A_0)\nabla\bar y(x'),\nabla\bar y(x')\ran\les
\lan(A_1-A_0)\nabla\bar y(x),\nabla\bar y(x)\ran\les\lan(A_1-A_0)\nabla\bar y(x''),\nabla y(x'')\ran,\\
\ns\ds\qq\qq\qq\ae x'\in\big(\bar\si(\cd)=0\big),\q\ae
x\in\big(0<\bar\si(\cd)<1\big),\q\ae x''\in\big(\bar\si(\cd)=1\big).\ea\right.\ee
Further, in the case
\bel{<b}\int_\O\bar\si(x)dx<\b|\O|,\ee
the following holds:
\bel{<0}\lan(A_1-A_0)\nabla\bar y(x),\nabla\bar y(x)\ran\les0,\qq\ae
x\in\big(0\les\bar\si(\cd)<1\big);\ee
in the case
\bel{>a}\int_\O\bar\si(x)dx>\a|\O|,\ee
the following holds:
\bel{>0}\lan(A_1-A_0)\nabla\bar y(x),\nabla\bar y(x)\ran\ges0,\qq\ae
x\in\big(0<\bar\si(\cd)\les1\big);\ee
and in the case
\bel{a<b}\a|\O|<\int_\O\bar\si(x)dx<\b|\O|,\ee
the following holds:
\bel{=0}\lan(A_1-A_0)\nabla\bar y(x),\nabla\bar y(x)\ran=0,\qq\ae
x\in\big(0<\bar\si(\cd)<1\big).\ee
\et

\it Proof. \rm Since $A\big(\Si[\a,\b]\big)$ is convex and closed in $L^1(\O;\dbS^n)$, by Corollary \ref{Corollary 2.7} part (i), we have that Problem ($\bar\L^c[\a,\b]$) admits an optimal solution $\bar\si(\cd)$ and the following is its characterization:
\bel{3.12}0\ges\int_\O\big(\si(x)-\bar\si(x)\big)\lan(A_1-A_0)\nabla
\bar y(x),\nabla\bar y(x)\ran dx,\qq\forall\si(\cd)\in\Si[\a,\b].\ee
We now look at further necessary conditions for $(\bar y(\cd),\bar\si(\cd))$.

\ms

Since \eqref{(0<si<1)>0} holds, for any $v(\cd)\in\sV_0$ where
$$\sV_0\equiv\Big\{v(\cd)\in L^\infty(\O;\dbR)\bigm|\int_\O v(x)dx=0\Big\},$$
with
$$\supp v(\cd)\subseteq(0<\bar\si(\cd)<1),$$
one has $\si(\cd)=\bar\si(\cd)\pm\e v(\cd)\in\Si[\a,\b]$ as long as
$\e>0$ is small enough. By taking such a $\si(\cd)$ in
(\ref{3.12}), we have
\bel{3.13}\int_\O v(x)\lan(A_1-A_0)\nabla\bar y(x),\nabla\bar y(x)\ran
dx=0,\qq\forall v(\cd)\in\sV_0.\ee
This leads to the first identity in (\ref{3.3A}).

\ms

Next, let the set $(\bar\si(\cd)=0)$ have a positive measure. Then take any $v(\cd)\in L^\infty(\O)$ with $v(x)\ges0$, supported on $(\bar\si(\cd)=0)$, and any $w(\cd)\in L^\infty(\O)$ with $w(x)\ges0$,
supported on $(0<\bar\si(\cd)\les1)$ (which has a positive measure by \eqref{(0<si<1)>0}), and
\bel{intw=intv}\int_{(0<\bar\si(\cd)\les1)}w(x)dx=\int_{(\bar\si(\cd)=0)}v(x)dx.\ee
Then for $\e>0$ small enough,
$$\si(\cd)=\bar\si(\cd)+\e v(\cd)-\e w(\cd)\in\Si[\a,\b].$$
Hence, using such a $\si(\cd)$ in (\ref{3.12}), one obtains the first inequality in the second conclusion of \eqref{3.3A}. Likewise, we can obtain the second inequality in the second conclusion of \eqref{3.3A}.

\ms

Further, if \eqref{<b} holds, we may take $v(\cd)\in L^\infty(\O)$ with $v(\cd)\ges0$, supported in $(0\les\bar\si(\cd)<1)$. Then, for $\e>0$ small, $\bar\si(\cd)+\e v(\cd)\in\Si[\a,\b]$. Taking such a $\si(\cd)$ in \eqref{3.12}, we obtain \eqref{<0}. Likewise we can obtain \eqref{>0} under \eqref{>a}. Finally, combining the above two cases, we obtain \eqref{=0} under \eqref{a<b}. This completes the proof. \endpf

We now present an interesting corollary.

\bc{A_0<A_1} \sl Let $0<\a<\b<1$, and $A_0<A_1$. Then there exists an optimal control $\bar\si(\cd)$ of Problem {\rm($\bar\L^c[\a,\b]$)} such that
\bel{int=b}\int_\O\bar\si(x)dx=\b|\O|.\ee
Further, if  $\bar\si(\cd)$ is a constant on $\O$, then $\bar\si(\cd)$ is not an optimal control of Problem {\rm($\bar\L^c[\a,\b]$)}.

\ec

\it Proof. \rm
 The first conclusion is obvious. It follows easily from the monotonicity of the principal eigen map (see \eqref{a<bar a}).

Next, suppose $\bar\si(x)\equiv\si_0\in[\a,\b]\subset(0,1)$ (so that $\bar\si(\cd)\in\Si[\a,\b]$). If such a $\bar\si(\cd)$ is optimal, then by the first equation in \eqref{3.3A}, we have
\bel{=C}\lan(A_1-A_0)\nabla\bar y(x),\nabla\bar y(x)\ran\equiv C,\qq\ae~x\in\O.\ee
On the other hand, since $\ds A(\bar\si(x))$ is a constant matrix, $\bar y(\cd)\in C^\infty(\O)$. Thus making use of the fact that $\bar y(\cd)$ is strictly positive in $\O$ and is zero on $\pa\O$, one sees that $\bar y(\cd)$ attains its maximum on $\ol\O$ at some point $x_0\in\O$. This implies $\na \bar y(x_0)=0$ and therefore,
$$C=\lan(A_1-A_0)\nabla\bar y(x_0),\nabla\bar y(x_0)\ran=0.$$
Since $A_1-A_0$ is positive definite, $C=0$ and \eqref{=C} implies $\bar y(\cd)\equiv 0$. This is a contradiction since $\bar y(\cd)$ is positive in $\O$. Hence, such a $\bar\si(\cd)$ is not optimal. \endpf

Note that a constant $\bar\si(\cd)\equiv\si_0$ is not an optimal control means that a perfect mixture of two different material does not gives the optimal solution to the problem.

\ms

To conclude this section, we present an illustrative example for the case
$\a=0$, $\b=1$ with both $A_0\les A_1$ and $A_0\ges A_1$ fail.

\bex{Ex3.3.} \rm Let $\O=[-1,1]\times[-1,1]$ which is a bounded Lipschitz domain, $\a=0$, $\b=1$, and
$$A_0=\begin{pmatrix}{1\over2}&0\\ 0&{3\over2}\end{pmatrix},\qq A_1=\begin{pmatrix}1&0\\0&1\end{pmatrix}.$$
Thus, both $A_0\les A_1$ and $A_0\ges A_1$ fail. Let us consider the following problem:
$$\left\{\2n\ba{ll}
\ds-\big(ay_{x_1x_1}+by_{x_2x_2}\big)=\l y,\qq\hb{in }~\O,\\
\ns\ds y\big|_{\pa\O}=0,\ea\right.$$
for any given $a,b>0$. Then we can check directly that the principal eigen pair is given by
$$\left\{\ba{ll}
\ns\ds\l={(a+b)\pi^2\over4},\\
\ns\ds y(x_1,x_2)=\cos{\pi x_1\over2}\cos{\pi x_2\over2},\qq(x_1,x_2)\in[-1,1]\times[-1,1].\ea\right.$$
From this, we see that for any constant $\si\in[0,1]$, one has
$$A_0+\si(A_1-A_0)=\begin{pmatrix}{1\over2}+{\si\over2}&0\\ 0&{3\over2}-{\si\over2}\end{pmatrix},$$
and by the above calculation,
\bel{l=l}\l_{A_0+\si(A_1-A_0)}={\pi^2\over2}=\l_{A_0}=\l_{A_1},\qq\forall\si
\in[0,1].\ee
This shows that the map $a(\cd)\mapsto\l_{a(\cd)}$ is not strictly convex. If $\si(x)\equiv\si\in(0,1)$ is optimal, then by Theorem \ref{Theorem 3.1}, we should have
$$\ba{ll}
\ns\ds C\equiv\lan(A_1-A_0)\nabla y(x),\nabla y(x)\ran
={\pi^2\over4}\lan\begin{pmatrix}{1\over2}&0\\ 0&-{1\over2}\end{pmatrix}\begin{pmatrix}\sin{\pi x_1\over2}\cos{\pi x_2\over2}\\ \cos{\pi x_1\over2}\sin{\pi x_2\over2}\end{pmatrix},\begin{pmatrix}\sin{\pi x_1\over2}\cos{\pi x_2\over2}\\ \cos{\pi x_1\over2}\sin{\pi x_2\over2}\end{pmatrix}\ran\\ [4mm]
\ns\ds={\pi^2\over8}\[\sin^2{\pi x_1\over2}\cos^2{\pi x_2\over2}-\cos^2{\pi x_1\over2}\sin^2{\pi x_2\over2}\]={\pi^2\over8}\sin{\pi(x_1+x_2)\over2}\sin{\pi(x_1-x_2)\over2},\ea$$
which is impossible. Hence, $\si(x)\equiv\si\in(0,1)$ is not an optimal solution to the corresponding maximization problem. Due to \eqref{l=l}, we see that both $\si(\cd)=0$ and $\si(\cd)=1$ are not optimal either. By the way, the above also roughly means that in the current case, if $A_0$ and $A_1$ represent the heat diffusibility of the two material, then the uniform mixture of any ratio of these two material is not optimal for Problem ($\bar\L^c[0,1]$). It is not clear to us at the moment what is an optimal control for this problem.

\ex

\section{Relaxation of Problem $(\underline\L[\a,\b])$}

Mimicking Problem ($\bar\L^c[\a,\b]$), we may pose the following problem.

\ms

\bf Problem ($\underline\L^c[\a,\b]$). \rm Let $0<\m_0\les\m_1<\infty$, $A_0,A_1\in M[\m_0,\m_1]$ and $0\les\a\les\b\les1$ satisfying (\ref{ab}). Find a $\underline\si(\cd)\in\Si[\a,\b]$ such that
\bel{}\l_{\underline\si(\cd)}=\inf_{\si(\cd)\in\Si[\a,\b]}\l_{\si(\cd)}.\ee

\ms

Note that although $\Si[\a,\b]$ is convex and closed (in $L^1(\O;\dbR)$), $A(\cd)\mapsto\l_{A(\cd)}$ is concave and not necessarily convex in general. Therefore, it is not clear if the map $\si(\cd)\mapsto\l_{\si(\cd)}$ admits a minimum on $\Si[\a,\b]$. In another word, the above Problem ($\underline\L^c[\a,\b]$) might not admit a minimum over $\Si[\a,\b]$ in general. Hence, instead of Problem ($\underline\L^c[\a,\b]$), we will introduce another relaxation of Problem ($\underline\L[\a,\b]$), for which the existence of
an optimal solution is guaranteed. To this end, let us recall some results relevant to the so-called $H$-convergence, which will play an essential role in the relaxation of Problem ($\underline\L[\a,\b]$).

\subsection{$H$-Convergence}

We recall the following definition.

\bde{H-convergence} \rm A sequence $\{a_\e(\cd)\}\subseteq\sM[\m_0,\m_1]$ is said to be $H$-{\it convergent} to $a^*(\cd)\in\sM[\m_0,\m_1]$ on $\O$, denoted by $a_\e(\cd)\hcon a^*(\cd)$, if for any $f\in W^{-1,2}(\O)$, the weak solution $y_\e(\cd)$ of the following problem
\bel{E401}\left\{\2n\ba{ll}
\ds-\nabla\cd\big(a_\e(x)\nabla y_\e(x)\big)=f,\qq\hb{in }\O,\\
\ns\ds y_\e\big|_{\pa\O}=0 \ea\right.\ee
has the property that
\bel{4.2}y_\e(\cd)\to y^*(\cd),\qq\hb{weakly in }\, W^{1,2}_0(\O)\ee
with $y^*(\cd)$ being the weak solution to the following:
\bel{4.3}\left\{\2n\ba{ll}
\ds-\nabla\cd\big(a^*(x)\nabla y^*(x)\big)=f,\qq\hb{in }\O,\\
\ns\ds y^*\big|_{\pa\O}=0. \ea\right.\ee

\ede

In 1968, Spagnolo (\cite{Spagnolo-1968}) introduced the above notion, called the {\it G-convergence}, for symmetric operators (i.e., each $a_\e(\cd)$ is  symmetric matrix valued and so is $a^*(\cd)$). The notion was generalized by Tartar for possibly non-symmetric operators (\cite{Tartar-1978}), and is called the $H$-convergence, for which the following additional condition is required:
\bel{4.2*}a_\e(\cd)\nabla y_\e(\cd)\to a^*(\cd)\nabla y^*(\cd),\qq
\hb{weakly in }\, L^2(\O;\dbR^n),\ee
which is automatically true when $a_\e(\cd)$ is symmetric and is $H$-convergent. It is known that for symmetric operators, the $G$-convergence is equivalent to the $H$-convergence (\cite{Allaire-2002}). In the problems that we are studying, all the involved second order differential operators are symmetric. Hence, $G$-convergence will be enough. However, we prefer to use the name $H$-convergence instead, just keep in mind that we are treating the case of symmetric operators.

\ms

Note that in the definition, the $H$-limit $a^*(\cd)$ of $a_\e(\cd)$ is independent of the choice of $f\in W^{-1,2}(\O)$, and the whole sequence (not just a subsequence) $y_\e(\cd)$ is required to be weakly convergent in $W^{1,2}_0(\O)$.

\ms

For any $X\subseteq\sM[\m_0,\m_1]$, we denote the $H$-{\it closure} of $X$ by $\cl X^H$. Let us now list some useful properties of $H$-convergence and $H$-closure, found in \cite{Allaire-2002}.

\ms

(i) {\bf Sequential compactness.} The set $\sM[\m_0,\m_1]$ is sequentially compact under $H$-convergence, i.e., for any sequence $\{a_k(\cd)\}_{k\ges1}\subseteq\sM[\m_0,\m_1]$, there exists a subsequence $\{a_{k_i}(\cd)\}_{i\ges1}$ and an $a^*(\cd)\in\sM[\m_0,\m_1]$ such that $a_{k_i}(\cd)\hcon a^*(\cd)$.

\ms

(ii) {\bf Locality.} If $a_\e(\cd)\hcon a^*(\cd)$ on $\O$, then $a_\e(\cd)\hcon a^*(\cd)$ on any subdomain $\o\subseteq\O$.

\ms

(iii) {\bf Monotonicity and uniqueness.} If $\ds a_\e(\cd)\hcon a^*(\cd)$, $\ds b_\e(\cd)\hcon b^*(\cd)$ with $\ds a_\e(\cd)\les b_\e(\cd)$, then $\ds a^*(\cd)\les b^*(\cd)$. In particular, if $0<\m_0\les a(x)\les\m_1$, then
$$\m_0\les a^*(x)\les\m_1,\qq x\in\O.$$
Also, by taking $b_\e(\cd)=a_\e(\cd)$, one has
$$\left.\ba{l}
\ds a_\e(\cd)\hcon a^*(\cd)\\
\ns\ds a_\e(\cd)\hcon b^*(\cd)\ea\right\}\q\Ra\q a^*(\cd)=b^*(\cd).$$
That is, the $H$-limit of a sequence is unique.

\ms

(iv) {\bf Non-homogeneous boundary conditions.} Let
$$a_\e(\cd)\hcon a^*(\cd),\q\f\in W^{1,2}(\O),\q f\in W^{-1,2}(\O),$$
and $y_\e(\cd)$ and $y^*(\cd)$ respectively be the solutions to the following:
\bel{E206}\left\{\2n\ba{ll}
\ds-\nabla\cd\big(a_\e(x)\nabla y_\e(x)\big)=f,\qq\hb{in }\O,\\
\ns\ds y_\e\big|_{\pa\O}=\f,\ea\right.\ee
\bel{E207}\left\{\2n\ba{ll}
\ds-\nabla\cd\big(a^*(x)\nabla y^*(x)\big)=f,\qq\hb{in }\O,\\
\ns\ds y^*\big|_{\pa\O}=\f.\ea\right.\ee
Then
$$y_\e(\cd)\to y^*(\cd),\qq\hb{weakly in }W^{1,2}(\O).$$

(v) {\bf Metrizability.} Let $\{\Bf_\ell(\cd)\}_{\ell\ges1}\subseteq L^2(\O;\dbR^n)$ such that $\{\nabla\cd\Bf_\ell\}_{\ell\ges1}$ is dense in $W^{-1,2}(\O)$. For any $a(\cd),b(\cd)\in\sM[\m_0,\m_1]$, let $y_\ell^{a(\cd)}(\cd)$ and $y^{b(\cd)}_\ell(\cd)$ be the unique weak solutions to the following:
$$\left\{\ba{ll}
\ds-\nabla\cd\(a(x)\nabla y^{a(\cd)}_\ell(x)\)=\nabla\cd\Bf_\ell,\qq\hb{in }\O,\\
\ns\ds y^{a(\cd)}_\ell\big|_{\pa\O}=0,\ea\right.$$
and
$$\left\{\ba{ll}
\ds-\nabla\cd\(b(x)\nabla y^{b(\cd)}_\ell(x)\)=\nabla\cd\Bf_\ell,\qq\hb{in }\O,\\
\ns\ds y^{b(\cd)}_\ell\big|_{\pa\O}=0.\ea\right.$$
Define
\bel{rho}\rho\big(a(\cd),b(\cd)\big)=\sum_{\ell\ges1}2^{-\ell}
{\|y^{a(\cd)}_\ell(\cd)
-y^{b(\cd)}_\ell(\cd)\|_2+\|a(\cd)\nabla y^{a(\cd)}_\ell(\cd)-b(\cd)\nabla y^{b(\cd)}_\ell(\cd)\|_2\over\|\Bf_\ell\|_2}.\ee
Then $\rho(\cd\,,\cd)$ is a metric on $\sM[\m_0,\m_1]$ such that for any $a_\e(\cd),a^*(\cd)\in\sM[\m_0,\m_1]$,
\bel{rho(a_e,a)}a_\e(\cd)\hcon a^*(\cd)\qq\iff\qq\rho\big(a_\e(\cd),a^*(\cd)\big)\to0.\ee

\ms

(vi) {\bf Upper and lower bounds.} Let
$$\left\{\2n\ba{ll}
\ds a_\e(\cd)\to\bar a(\cd),\q a_\e(\cd)^{-1}\to\underline a(\cd)^{-1},\q\;\hb{weak$^*$ in }L^\infty(\O;\dbS^n),\\
\ns\ds a_\e(\cd)\hcon a^*(\cd).\ea\right.$$
Then
\bel{E205}\underline a(x)\les a^*(x)\les\bar a(x),\qq\ae x\in\O.\ee

\ms

(vii) {\bf Commutativity with congruent transformation.} Let $Q\in\dbR^{n\times n}$ be non-singular. Then
$$a_\e(\cd)\hcon a^*(\cd)\q\iff\q Q a_\e(\cd)Q^\RT\hcon Qa^*(\cd)Q^\RT.$$

\ms

(viii) {\bf Pointwiseness.} Let $G\subseteq M[\m_0,\m_1]$, denote
$$L^\infty(\O;G)=\Big\{a(\cd)\in L^\infty(\O;\dbS^n)\bigm|a(x)\in G,~\ae x\in\O\Big\},$$
and define
$$\cl G^{\,H}=\Big\{A\in M[\m_0,\m_1]\bigm|\chi_{_\O}(\cd)A\in \cl{L^\infty(\O;G)}^{\,H}\Big\}.$$
Then
\bel{L=L}\cl{L^\infty(\O;G)}^{\,H}=L^\infty\big(\O;\cl G^H\big)\equiv\Big\{a^*(\cd)\in\sM[\m_0,\m_1]\bigm|a^*(x)\in\cl G^{\,H},~\ae x\in\O\Big\}.\ee
Namely, $a^*(\cd)\in\cl{L^\infty(\O;G)}^{\,H}$ if and only if for almost all $x\in\O$, there exists a sequence $\{a_k(\cd\,;x)\}_{k\ges1}\subseteq L^\infty(\O;G)$ (depending on $x$) such that
$$a_k(\cd\,;x)\hcon \chi_\O(\cd)a^*(x),\qq k\to\infty.$$
More generally, let $\cQ\subseteq\sM[\m_0,\m_0]$ and define
$$Q_x=\big\{q(x)\bigm|q(\cd)\in\cQ\big\},\qq\forall x\in\O.$$
Then, under some mild conditions (see Theorem 2.3 in \cite{Li-Lou-2011})
$$\cl\cQ^{\,H}=\Big\{a^*(\cd)\in\sM[\m_0,\m_1]\bigm|a^*(x)\in\cl Q_x^{\,H},~\ae x\in\O\Big\}.$$
Note that by taking $a_k(\cd)=\chi_{_\O}(\cd)A$ with $A\in G$, we see that
\bel{G in G}G\subseteq\cl G^H.\ee
We will see that $G$ is a proper subset of $\cl G^H$ below.

\subsection{Lamination}

In this subsection, we consider a special case involving two matrices, which will be useful in our relaxation of Problem ($\underline\L[\a,\b]$). Let us first present the following result.

\bt{Theorem 4.2} \sl Let $0<\m_0\les\m_1<\infty$ and $A,B\in M[\m_0,\m_1]$ be fixed.

\ms

{\rm(i)} For any $\th\in(0,1)$ and $e\in S^{n-1}\equiv\{x\in\dbR^n\bigm||x|=1\}$, define
\bel{H_e}\cH_\e(A,B;\th,e)=\left\{\2n\ba{ll}
\ds A,\qq\Big\{{\lan x,e\ran\over\e}\Big\}\in[\th,1),\\
\ns\ds B,\qq\Big\{{\lan x,e\ran\over\e}\Big\}\in[0,\th),\ea\right.\qq x\in\dbR^n,\ee
where $\{r\}=r-[r]$ is the decimal part of the real number $r$. Then, as $\e\to 0^+$,
\bel{A*}\ba{ll}
\ns\ds\cH_\e(A,B;\th,e)\hcon\cH[A,B;\th,e]\\
\ns\ds\qq\qq\equiv(1-\th)A+\th B-{\th(1-\th)(A-B)ee^\RT
(A-B)\over e^\RT\big[\th A+(1-\th)B\big]e}\\
\ns\ds\qq\qq=A-\th(A-B)-{\th(1-\th)(A-B)ee^\RT
(A-B)\over e^\RT\big[B+\th(A-B)\big]e}\in\cl{\{A,B\}}^{\,H}\subseteq M[\m_0,\m_1].\ea\ee

{\rm(ii)} For any $m\ges1$, let
\bel{Gm}\ba{ll}
\ns\ds\G^m(A,B)=\Big\{A^*\in\dbS^n\bigm|(1-\th)(A-B)=(A^*-B)\[I+\th\sum_{k=1}^m\b_k{e_ke_k^\RT
(A-B)\over e_k^\RT Be_k}\],\\
\ns\ds\qq\qq\qq\qq\qq\qq\hb{ for some }\th\in[0,1],~e_k\in S^{n-1},~\b_k\ges0,~\sum_{k=1}^m\b_k=1\Big\},\ea\ee
and
\bel{G(A,B)}\G(A,B)=\bigcup_{m=1}^\infty\G^m(A,B)\subseteq\cl{\{A,B\}}^{\,H}.\ee
Then
\bel{G1=}\G^1(A,B)=\Big\{\cH[A,B;\th,e]\bigm|\th\in[0,1],~e\in S^{n-1}\Big\}\equiv\cH\big[A,B;[0,1],
S^{n-1}\big],\ee
and
\bel{G(G)inG}\G^1\big(\G(A,B),B\big)\subseteq\G(A,B).\ee

{\rm(iii)} For any $\th\in[0,1]$, $H\in\dbS^n$ with $H\ges0$ and $\tr(H)=1$, the matrix
$$I+\th B^{-{1\over2}}HB^{-{1\over2}}(A-B)$$
is non-singular, and
\bel{E409}\ba{ll}
\ns\ds\G(A,B)=\Big\{B+(1-\th)(A-B)\big[I+\th B^{-{1\over2}}H
B^{-{1\over2}}(A-B)\big]^{-1}\bigm|\\
\ns\ds\qq\qq\qq\qq\qq\qq\qq\qq\qq\th\in[0,1],~H\ges0,~\tr(H)=1\Big\}.\ea\ee

\et

\it Proof. \rm (i) Relation \eqref{A*} follows from \cite{Allaire-2002}, Corollary 1.3.34 (see \cite{Lou-Yong-2009} also).

\ms

(ii) This is a restatement of Lemma 2.2.3 of \cite{Allaire-2002}.

\ms

(iii) Let $\th\in(0,1)$ and $H>0$, $\tr(H)=1$, we have $I-H\ges0$, and
%
$$\ba{ll}
\ns\ds\det\Big[I+\th B^{-{1\over2}}H B^{-{1\over2}}(A-B)\Big]\\
\ns\ds=\det(B^{1\over2})\det\Big[I+\th B^{-{1\over2}}H B^{-{1\over2}}A-\th B^{-{1\over2}}HB^{1\over2}\Big]\det(B^{-{1\over2}})\\
\ns\ds=\det\Big[I+\th H( B^{-{1\over2}}AB^{-{1\over2}}-I)\Big]\\
\ns\ds=\det(H^{-{1\over2}})\det\Big[I-H+(1-\th)H+\th HB^{-{1\over2}}AB^{-{1\over2}}\Big]
\det(H^{1\over2})\\
\ns\ds=\det\Big[I-H+(1-\th)H+\th H^{1\over 2}B^{-{1\over2}}AB^{-{1\over2}}H^{1\over 2}\Big]\\
\ns\ds\ges\det\Big[I-H+(1-\th)H+\th{\m_0\over\m_1}H\Big]=\det\[I-\th{\m_1-\m_0\over\m_1}H\]\\
\ns\ds=\det\[{\th\m_0+(1-\th)\m_1\over\m_1}I+{\th(\m_1-\m_0)\over\m_1}(I-H)\]\ges
\[{\th\m_0+(1-\th)\m_1\over\m_1}\]^n.\ea$$
By continuity, for any $\th\in[0,1]$, $H\ges0$ (instead of just $H>0$), with $\tr H=1$, we have
$$\det\1n\Big[I+\th B^{-{1\over2}}H B^{-{1\over2}}(A\1n-\1n B)\Big]
\ges
\[{\th\m_0+(1-\th)\m_1\over\m_1}\]^n>0.$$
Therefore, $\Big[I+\th B^{-{1\over2}}HB^{-{1\over2}}(A-B)\Big]$ is non-singular.

\ms

On the other hand, if $A^*\in\G(A,B)$, then
\bel{4.22}(1-\th)(A-B)=(A^*-B)\[I+\th\sum^m_{k=1}\b_k{e_k e_k^\RT(A-B)\over e_k^\RT B e_k}\],\ee
for some $\th\in[0,1]$, $\b_k\ges0$, $\ds\sum_{k=1}^m\b_k=1$, and $e_k\in S^{n-1}$. Note that
$$\sum_{k=1}^m\b_k{e_ke_k^\RT\over e_k^\RT Be_k}=B^{-{1\over2}}\sum_{k=1}^m\b_k{(B^{1\over2}e_k)\over|B^{1\over2}e_k|}\({(B^{1\over2}e_k)
\over|B^{1\over2}e_k|}\)^\RT B^{-{1\over2}}\equiv B^{-{1\over2}}HB^{-{1\over2}},$$
where
$$H=\sum_{k=1}^m\b_k{(B^{1\over2}e_k)^\RT\over|B^{1\over2}e_k|}\({(B^{1\over2}e_k)
\over|B^{1\over2}e_k|}\)^\RT\ges0,\qq\tr(H)=\sum_{k=1}^m\b_k=1.$$
Thus, \eqref{4.22} is equivalent to the following:
\bel{4.23}(1-\th)(A-B)=(A^*-B)\[I+\th B^{-{1\over2}}HB^{-{1\over2}}(A-B)\].\ee
Then by the invertibility of $\[I+\th B^{-{1\over2}}HB^{-{1\over2}}(A-B)\]$, we have
\bel{4.24}A^*=B+(1-\th)(A-B)\[I+\th B^{-{1\over2}}HB^{-{1\over2}}(A-B)\]^{-1}.\ee
Conversely, if \eqref{4.24} holds for some $\th\in[0,1]$, and $H\ges0$, $\tr(H)=1$, then \eqref{4.23} holds. Moreover, it is easy to see that
$$H=\sum^n_{k=1}\b_k \xi_k\xi_k^\RT$$
with $\b_k\ges0$ $(1\les k\les n$), $\ds\sum^n_{k=1}\b_k=1$ and $\xi_k\in S^{n-1}$  $(1\les k\les n$). Thus
$$B^{-{1\over2}}HB^{-{1\over2}}=\sum^n_{k=1}\b_k{\xi_k\xi_k^\RT\over e_k^\RT B e_k},$$
where
$$e_k={B^{-{1\over2}}\xi_k\over|B^{-{1\over2}}\xi_k|}\in S^{n-1},\q 1\les k\les n.$$
Hence,
$$\ba{ll}
\ns\ds(1-\th)(A-B)=(A^*-B)\[I+\th B^{-{1\over2}}\sum_{k=1}^n\b_k{\xi_k\xi_k^\RT\over e_k^\RT Be_k}B^{-{1\over2}}(A-B)\]\\
\ns\ds\qq\qq\qq\q=(A^*-B)\[I+\th\sum_{k=1}^n\b_k{e_ke_k^\RT\over e_k^\RT Be_k}(A-B)\],\ea$$
which means $A^*\in\G(A,B)$. This completes the proof. \endpf

In the above, any element in $\G(A,B)$ is called a {\it lamination} of $A$ with base $B$. From \eqref{E409}, we see that $A,B\in\G(A,B)$ (by taking $\th=0,1$). Thus,
$$\{A,B\}\subsetneq\G(A,B)\subseteq\cl{\{A,B\}}^H.$$
Note also that for any $\th\in(0,1)$, $H\ges0$, $\tr(H)=1$ it holds that (see (\ref{E205}))
\bel{E208}\ba{ll}
\ns\ds\Big((1-\th)A^{-1}+\th B^{-1}\Big)^{-1}\les B+(1-\th)(A-B)\Big[I+\th B^{-{1\over 2}}H B^{-{1\over
2}}(A-B)\Big]^{-1}\les(1-\th)A+\th B.\ea\ee
This gives bounds for elements in $\G(A,B)\subseteq M[\m_0,\m_1]$. Further, we should keep in mind some facts about the set $\G(A,B)$:

\ms

$\bullet$ It is possible that $\G(A,B)\ne\G(B,A)$ ($n\ges 3$)

\ms

$\bullet$ $\G(A,B)$ is not necessarily convex, and even $\cl{\{A,B\}}^{\,H}$ might be non-convex.

\ms

$\bullet$ It is possible that $\cl{\{A,B\}}^{\,H}\ne\G(A,B)\bigcup\G(B,A)$.

\ms

$\bullet$ Even for $A=\l B$ with $B$ being diagonal ($n\ges2$ and $B\ne\g I$ for any $\g\in\dbR$, of course), as long as $\l\ne1$, $\G(A,B)$ contains non-diagonal matrices.

\subsection{Relaxation problem}

In this subsection, we fix $0<\m_0\les\m_1<\infty$, $0\les\a\les\b\les1$, and $G=\{A_0,A_1\}\subseteq M[\m_0,\m_1]$ satisfying \eqref{ab}.
 For a domain $\o$, denote
$$
\sA_{\o}[\a,\b]=\1n\Big\{A_0+\chi_{_{\O_1}}(\cd)(A_1-A_0)
\bigm|\O_1\subseteq \o \hb{ measurable, }~\a|\o|\1n\les\1n|\O_1|\1n\les\1n\b|\o|\Big\}.$$
Recall (see \eqref{sA[a,b]}) that
$\ds \sA[\a,\b]= \sA_{\O}[\a,\b]$.
We first present a simple result.

\bp{simple} \sl Assume $A_0\ne A_1$.  The following hold:
\bel{A_0notin}\a=0\qq\iff\qq\chi_{_\O}(\cd)A_0\in\cl{\sA[\a,\b]}^{\,H},\ee
\bel{A_1notin}\b=1\qq\iff\qq\chi_{_\O}(\cd)A_1\in\cl{\sA[\a,\b]}^{\,H}.\ee

\ep

\it Proof. \rm We just prove \eqref{A_0notin}. The other is similar.

\ms

Suppose $\a=0$. Then
$$\chi_{_\O}(\cd)A_0\in\sA[0,\b]\subseteq\cl{\sA[0,\b]}^{\,H}.$$
Conversely, suppose $\chi_{_\O}(\cd)A_0\in\cl{\sA[\a,\b]}^{\,H}$. Then there exists a sequence
$$a_k(\cd)=\chi_{_{\O_k^c}}(\cd)A_0+\chi_{_{\O_k}}(\cd)A_1\in\sA[\a,\b],\qq
\a|\O|\les|\O_k|\les\b|\O|,\qq k\ges1,$$
such that $a_k(\cd)\hcon \chi_{\O}(\cd)A_0$. We may let
$$\chi_{_{\O_k}}(\cd)\to g(\cd),\qq\hb{weak$^*$ in }L^\infty(\O;\dbS^n),$$
with
$$s\equiv {1\over |\O|}\int_\O g(x)dx=\lim_{k\to\infty}{|\O_k|\over |\O|}\in[\a,\b].$$
Then by (vi) of listed properties of $H$-convergence, we have
$$A_0\les[1-g(\cd)]A_0+g(\cd)A_1,\qq A_0^{-1}\les [1-g(\cd)]A_0^{-1}+g(\cd)A_1^{-1}.$$
Thus, integrating each side, one has
$$A_0\les(1-s)A_0+sA_1,\qq A_0^{-1}\les(1-s)A_0^{-1}+sA_1^{-1}.$$
Since $A_0\ne A_1$, the above hold only if $s=0$. Hence, $\a=0$, proving the conclusion. \endpf

\ms

We now formulate the following problem which is called an $H$-relaxation of Problem ($\underline\L[\a,\b]$), with $H$ indicating that the relaxation is in the sense of $H$-convergence.

\ms

\bf Problem ($\underline\L^H[\a,\b]$). \rm Find an $\underline a(\cd)\in\cl{\sA[\a,\b]}^{\,H}$ such that
\bel{E304}\l_{\underline a(\cd)}=\inf_{a(\cd)\in\cl{\sA[\a,\b]}^{\,H}}\l_{a(\cd)}\ee

Any $\underline a(\cd)$ satisfying \eqref{E304} is called an optimal control of Problem ($\underline\L^H[\a,\b]$), which is also called an optimal $H$-relaxed control of Problem ($\underline\L[\a,\b]$). The superscript ``$H$'' indicates the $H$-relaxation. We first have the following existence theorem.

\bt{Theorem 4.3} \sl Problem {\rm($\underline\L^H[\a,\b]$)} admits an optimal control $\underline a(\cd)\in\cl{\sA[\a,\b]}^{\,H}$.

\et

\it Proof. \rm Let $a_k(\cd)\in\cl{\sA[\a,\b]}^{\,H}$ be a minimizing sequence of Problem {\rm($\underline\L^H[\a,\b]$)} with $(\l_k,y_k(\cd))$ being the corresponding principal eigen-pair. Thus,
\bel{E309}\left\{\2n\ba{ll}
\ds-\nabla\cd\big(a_k(x)\nabla y_k(x)\big)=\l_ky_k(x),\qq x\in\O,\\
\ns\ds y_k|_{\pa\O}=0, \ea\right.\ee
\bel{E310}y_k(x)\ges0,\q x\in\O,\qq\int_\O|y_k(x)|^2\, dx=1,\ee
and
$$\lim_{k\to\infty}\l_k=\underline\l\equiv\inf_{a(\cd)\in
\cl{\sA[\a,\b]}^{\,H}}\l_{a(\cd)}.$$
Since $y_k(\cd)$ is uniformly bounded in $W^{1,2}_0(\O)$, and $\sM[\m_0,\m_1]$ is sequentially compact under $H$-convergence (Property (i) of $H$-convergence listed in Subsection 4.1), we may suppose that
\bel{E311}y_k(\cd)\to\underline y(\cd),\qq\hb{weakly in }W^{1,2}_0(\O),\,\hb{strongly in }L^2(\O),\ee
and
\bel{E312}a_k(\cd)\hcon\underline a(\cd),\ee
for some $\underline y(\cd)\in W^{1,2}_0(\O)$ and $\underline a(\cd)\in\cl{\sA[\a,\b]}^{\,H}$. Clearly,
\bel{E313}\underline y(x)\ges0,\q x\in\O,\qq\int_\O|\underline y(x)|^2dx=1,\ee
and
\bel{E314}\left\{\2n\ba{ll}
\ds-\nabla\cd\big(\underline a(x)\nabla\underline y(x)\big)=\underline\l\,\underline y(x),\qq x\in\O,\\
\ns\ds\underline y\big|_{\pa\O}=0,\ea\right.\ee
Hence, $\l_{\underline a(\cd)}=\underline\l$ and $\underline a(\cd)\in\cl{\sA[\a,\b]}^{\,H}$ is an optimal control. \endpf

Now, we state the following necessary conditions for the optimal control of Problem ($\underline\L^H[\a,\b]$).

\bt{T601} \sl Let $\underline a(\cd)\in\cl{\sA[\a,\b]}^{\,H}$ be an optimal control of Problem {\rm($\underline\L^H[\a,\b]$)} with $(\underline\l,\underline y(\cd))\in[\m_0,\m_1]\times W^{1,2}_0(\O)$ being the corresponding principal eigen pair. Then
\bel{sup**}\ba{ll}
\ns\ds\underline\l=\int_\O|\underline a(x)^{1\over2}\nabla\underline y(x)|^2dx=\sup_{b(\cd)\in\cl{\sA[\a,\b]}^{\,H}}\int_\O|b(x)^{-{1\over2}}\underline a(x)\nabla\underline y(x)|^2dx\\
\ns\ds\qq\qq\qq\qq\qq\q=\sup_{b(\cd)\in\sA[\a,\b]}\int_\O|b(x)^{-{1\over2}}\underline a(x)\nabla\underline y(x)|^2dx,\ea\ee
and
\bel{EN601}\int_\O\ip{\big[\underline a(x)-\underline a(x)b(x)^{-1}\underline a(x)\big]\nabla\underline y(x),\nabla\underline y(x)}\,dx\ges0,\qq\all b(\cd)\in\cl{\sA[\a,\b]}^{\,H},\ee
Equivalently,
\bel{EN602}\int_\O\ip{\big[\underline a(x)-\underline a(x)b(x)^{-1}\underline a(x)\big]\nabla\underline y(x),\nabla\underline y(x)}\,dx\ges0,\qq\all b(\cd)\in\sA[\a,\b].\ee
When $\a=0$ and $\b=1$, the following also holds:
\bel{in A^H}\lan\underline a(x)\nabla\underline y(x),\nabla\underline y(x)\ran\ges\lan B^{-1}\underline a(x)\nabla\underline y(x),\underline a(x)\nabla\underline y(x)\ran,\q\ae x\in\O,~B\in\cl{\{A_0,A_1\}}^{\,H},\ee
which is equivalent to
\bel{in A}\lan\underline a(x)\nabla\underline y(x),\nabla\underline y(x)\ran\ges\lan A_i^{-1}\underline a(x)\nabla\underline y(x),\underline a(x)\nabla\underline y(x)\ran,\q\ae x\in\O,~i=0,1.\ee

\et

This is a kind of maximum principle for the optimal control of Problem ($\underline\L^H[\a,\b]$). We point out that when $0<\a$ or $\b<1$, one could not get \eqref{in A^H}. Also, it is clear that if \eqref{in A} holds, then by $H$-convergence,
\bel{in bar cA}\ba{ll}
\ns\ds\lan\underline a(x)\nabla\underline y(x),\nabla\underline y(x)\ran\ges\lan b(x)^{-1}\underline a(x)\nabla\underline y(x),\underline a(x)\nabla\underline y(x)\ran,\q\ae x\in\O,\q\forall b(\cd)\in\cl{\sA[\a,\b]}^{\,H}.\ea\ee
In other words, in some sense, \eqref{in A^H}, \eqref{in A} and \eqref{in bar cA} are mutually  equivalent.
To prove this theorem, we need several lemmas.

\bl{Lemma 4.6.} \sl The metric $\rho(\cd\,,\cd)$ on $\sM[\m_0,\m_1]$ defined by \eqref{rho} is uniformly continuous in the following sense: For any $\e>0$, there exists a $\d>0$ (only depending on $\e>0$) such that
\bel{rho<e}\rho(a(\cd),b(\cd))<\e,\qq\forall a(\cd),b(\cd)\in\sM[\m_0,\m_1],~\|a(\cd)-b(\cd)\|_1<\d.\ee
Consequently, if $\{a_k(\cd)\}$ and $\{b_k(\cd)\}$ are two sequences in $\sM[\m_0,\m_1]$ such that
\bel{ab_k}a_k(\cd)\hcon a^*(\cd),\qq\|a_k(\cd)-b_k(\cd)\|_1\to0.\ee
Then $b_k(\cd)\hcon a^*(\cd)$.

\el

\it Proof. \rm
Since $\sM[\m_0,\m_1]$ is sequentially compact under $H$-convergence, it
suffices to show that
\eqref{ab_k} implies $b_k(\cd)\hcon a^*(\cd)$.

\ms

Let $f\in W^{-1,2}(\O)$. Consider
\bel{eq-a}\left\{\2n\ba{ll}
\ds-\nabla\cd\big(a_k(x)\nabla y_k(x)\big)=f,\q\eqin\O,\\
\ns\ds y_k\big|_{\pa\O}=0\ea\right.\ee
and
\bel{eq-b}\left\{\2n\ba{ll}
\ds-\nabla\cd\big(b_k(x)\nabla z_k(x)\big)=f,\q\eqin\O,\\
\ns\ds z_k\big|_{\pa\O}=0.\ea\right.\ee
Then
$$y_k(\cd)\to y^*(\cd),\qq\wi{W^{1,2}_0(\O)}$$
with
$$\left\{\2n\ba{ll}
\ds-\nabla\cd\big(a^*(x)\nabla y^*(x)\big)=f,\q\eqin\O,\\
\ns\ds y^*\big|_{\pa\O}=0.\ea\right.$$
Then, thanks to Theorem \ref{Meyers-Gallouet-Monier}, for some $p>2$ and $C=C_f$, the following holds:
$$\|\na y_k\|_{L^p(\O)}\les C,\qq\|\na z_k\|_{L^p(\O)}\les C.$$
We note that (making use of the Dominated Convergence Theorem)
$$\|a_k(\cd)-b_k(\cd)\|_{2p\over p-2}\to0.$$
It holds that
\begin{eqnarray*}
0&=&\int_\O a_k(x)\na y_k(x)\cd  \na\big( y_k(x)-z_k(x)\big)\, dx- \int_\O  b_k(x)\na z_k(x)\cd\na\big( y_k(x)-z_k(x)\big)\, dx\\
&= & \int_\O a_k(x)\big(\na y_k(x)-\na z_k(x)\big)\cd \big(\na y_k(x)-\na z_k(x)\big)\, dx\\
&& +\int_\O \big(a_k(x)-b_k(x)\big)\na z_k(x)\cd\na\big( y_k(x)-z_k(x)\big)\Big)\, dx\\
&\ges & \m_0\|\na (y_k(\cd)-z_k(\cd)\|_2-\|a_k(\cd)-b_k(\cd)\|_{2p\over p-2}\|\na (z_k(\cd)\|_p\big(\|\na (y_k(\cd)-z_k(\cd))\|_p\big).
\end{eqnarray*}
Therefore,
$$
z_k(\cd)- y_k(\cd)\to 0,\qq\sti{W^{1,2}_0(\O)}.
$$
This implies
$$
z_k(\cd)\to y^*(\cd),\qq\wi{W^{1,2}_0(\O)}.
$$
That is, $b_k(\cd)\hcon a^*(\cd)$. This completes the proof. \endpf

\ms

The above result shows that there exists a non-decreasing function $h:[0,\infty)\to[0,\infty)$ with $h(0)=h(0+)=0$ such that
\bel{rho<h}\rho\big(a(\cd),b(\cd)\big)\les h\big(\|a(\cd)-b(\cd)\|_1\big),\qq\forall a(\cd),b(\cd)\in\sM[\m_0,\m_1].\ee
Such a relation will be used below.

\bl{Lemma 4.7N.}
{\rm(i)} Let $\set{\O_i}$ be a sequence of measurable subsets of $\O$ such that
$$
|\O_i|\to \g |\O|
$$
and
$$
a_i(\cd)=A_0+\chi_{_{\O_i}}(\cd)(A_1-A_0)\hcon a^*(\cd).
$$
Then $a^*(\cd)\in \cl{\sA[\g,\g]}^{\,H}$.

{\rm(ii)} For any $a^*(\cd)\in \cl{\sA[\a,\b]}^{\,H}$, there is a $\g\in [\a,\b]$ such that  $a^*(\cd)\in \cl{\sA[\g,\g]}^{\,H}$.
\el

\proof {\rm(i)} \rm For any $i$, we can choose a measurable set $\wt\O_i\subseteq\O$ such that $|\wt\O_i|=\g |\O|$ and
$$\int_\O\big|\chi_{\O_i}(x)-\chi_{\wt\O_i}(x)\big|dx=\Big||\O_i|-\g|\O|\Big|.$$
In fact, we can choose $\wt\O_i\subseteq\O_i$ if $|\O_i|\ges\g|\O|$ and $\wt \O_i\supset\O_i$ if $|\O_i|<\g|\O|$. Let
$$\wt a_i(\cd)=A_0+\chi_{_{\wt\O_i}}(\cd)(A_1-A_0)\in\sA[\g,\g].$$
Then $\|\wt a_i(\cd)-a_i(\cd)\|_1\to0$. By Lemma \ref{Lemma 4.6.}, we get $\wt a_i(\cd)\hcon a^*(\cd)$ and consequently  $a^*(\cd)\in\cl{\sA[\g,\g]}^{\,H}$.

\ms

{\rm(ii)} The result follows directly from  (i).
\endpf

\bl{Lemma 4.7.} \sl Let $\{\O_i\subseteq\O\bigm|1\les i\les m\}$ be a  partition of $\O$, i.e., it is a set of mutually disjoint domain such that
\bel{O=}\big|\O\setminus\bigcup_{i=1}^m\O_i\big|=0.\ee

{\rm(i)} For $i=1,2,\cds,m$, let $a_\e^i(\cd)\in\sM[\m_0,\m_1]$ be such that
\bel{E433}a_\e^i(\cd)\hcon a^i_0(\cd),\qq1\les i\les m.\ee
Then
\bel{E445}\sum_{i=1}^m \chi_{_{\O_i}}(\cd)a_\e^i(\cd)\hcon \sum_{i=1}^m\chi_{_{\O_i}}(\cd)a^i_0(\cd).\ee

{\rm(ii)} For $i=1,2,\cds,m$, let $\ds
a_i(\cd)\big|_{\O_i}\in\cl{\sA_{\O_i}[0,1]}^{\,H}$. Then
\bel{}\sum_{i=1}^m\chi_{_{\O_i}}(\cd)a_i(\cd)\in\cl{\sA[\a,\b]}^{\,H}\ee
if and only if for some $\g_1,\g_2,\ldots,\g_m\in[0,1]$, one has $\ds
a_i(\cd)\big|_{\O_i}\in\cl{\sA_{\O_i}[\g_i,\g_i]}^{\,H}$ and
\bel{}\a|\O|\les\sum^m_{i=1}\g_i|\O_i|\les\b|\O|.\ee

\el

\proof {\rm(i)} \rm
By the locality of $H$-convergence (see Property (ii) in \S 4.1), we get $a^i_\e(\cd)\hcon a^i_0(\cd)$  on $\O_i$ for every $i=1,2,\ldots m$.

\ms

On the other hand, by the compactness of $H$-convergence (see Property (i) in \S 4.1), along a subsequence, $\ds \sum_{i=1}^m\chi_{_{\O_i}}(\cd)a_\e^i(\cd)\hcon a^*(\cd)$  on $\O$ for some $a^*(\cd)\in \sM[\m_0,\m_1]$. Then by locality, $a^i_\e(\cd)\hcon a^*(\cd)$ on $\O_i$ for every $i=1,2,\ldots m$.
By the uniqueness of
$H$-convergent limit (see Property (iii) in \S 4.1), $a^*(\cd)=a^i_0(\cd)$ on  $\O_i$ ($1\les i\les m$). Thus, $\ds a^*
(\cd)=\sum_{i=1}^m\chi_{_{\O_i}}(\cd)a^i_0(\cd)$ on $\O$. Consequently, \eqref{E445} holds, not only in the sense of subsequence.

\ms

 {\rm(ii)}  \textbf{Sufficiency.} For  $i=1,2,\cds,m$, we have $\O_i^k\subset \O_i$ such that $|\O_i^k|=\g_i|\O_i|$ and
$$
a_i^k(\cd)=A_0+\chi_{_{\O_i^k}}(\cd)(A_1-A_0)\hcon a_i(\cd), \qq \eqon \O_i.
$$
Thus, by (i),
$$
 \sum_{i=1}^m \chi_{_{\O_i}}(\cd)a_i^k(\cd)\hcon \sum_{i=1}^m\chi_{_{\O_i}}(\cd)a_i(\cd).
$$
On the other hand,
$$
{1\over |\O|}\Big|\bigcup^m_{i=1}\O_i^k\Big|={1\over |\O|}\sum^m_{i=1}\g_i|\O_i|\equiv \g\in [\a,\b].
$$
Therefore, $\ds\sum_{i=1}^m\chi_{_{\O_i}}(\cd)a_i(\cd)\in \cl{\sA[\g,\g]}^{\,H}\subseteq \cl{\sA[\a,\b]}^{\,H}$.

\textbf{Necessity.} By Lemma \ref{Lemma 4.7N.}, there is $\g\in [\a,\b]$ such that  $\ds\sum_{i=1}^m\chi_{_{\O_i}}(\cd)a_i(\cd)\in \cl{\sA[\g,\g]}^{\,H}$.
Thus, there is  a sequence $\set{E_k}$  of measurable subsets of $\O$ such that
$\ds
|E_k|=\g |\O|
$
and
$$
a_k(\cd)=A_0+\chi_{_{E_k}}(\cd)(A_1-A_0)\hcon \sum_{i=1}^m\chi_{_{\O_i}}(\cd)a_i(\cd).
$$
Then, by locality, $a_k(\cd)\hcon a_i(\cd)$ on $\O_i$.  On the other had, we can suppose that for  every  $i
=1,2,\ldots,m$, $|E_k\cap \O_i|$ convergence to $\g_i|\O_i|$ for some $\g_i\in [0,1]$.
Thus  $\ds
a_i(\cd)\big|_{\O_i}\in\cl{\sA_{\O_i}[\g_i,\g_i]}^{\,H}$ and
$$
\sum^m_{i=1}\g_i|\O_i|=\lim_{k\to +\infty}\sum^m_{i=1} |E_k\cap \O_i|=\lim_{k\to +\infty} |E_k|=\g |\O|.
$$
We get the proof.
\endpf

\begin{remark} \rm We would like to mention that $\g$ in Lemma \ref{Lemma 4.7N.}{\rm (ii)} might be not unique. Therefore, it is possible that although $\ds\sum_{i=1}^m\chi_{_{\O_i}}(\cd)a_i(\cd)\in\cl{\sA[\a,\b]}^{\,H}$
and $\ds
a_i(\cd)\big|_{\O_i}\in\cl{\sA_{\O_i}[\g_i,\g_i]}^{\,H}$, but $\ds
{1\over|\O|}\sum^m_{i=1}\g_i|\O_i|\not\in[\a,\b]$.
\end{remark}

The following result is an extension of Theorem \ref{Theorem 4.2} (i), replacing $A$ and $B$ by $a(\cd)$ and $b(\cd)$, respectively.

\bl{Lemma 4.8.} \sl Let $a(\cd),b(\cd)\in\sM[\m_0,\m_1]$, $\th\in[0,1]$ and $e\in S^{n-1}$. For any small $\e>0$, define
\bel{}\cH_\e[a(\cd),b(\cd);\th,e](x)=\left\{\2n\ba{ll}
\ds a(x),\qq\Big\{{\lan x,e\ran\over\e}\Big\}\in[\th,1),\\
\ns\ds b(x),\qq\Big\{{\lan x,e\ran\over\e}\Big\}\in[0,\th).\ea\right.\ee
Then as $\e\to0$,
\bel{}\ba{ll}
\ns\ds\cH_\e[a(\cd),b(\cd);\th,e](\cd)\hcon\cH[a(\cd),b(\cd);\th,e](\cd)\\
\ns\ds\qq\qq\qq\qq\qq\equiv a(\cd)-\th[a(\cd)-b(\cd)]-{\th(1-\th)[a(\cd)-b(\cd)]ee^\RT[a(\cd)-b(\cd)]\over e^\RT\{b(\cd)+\th[a(\cd)-b(\cd)]\}e}.\ea\ee
\el

The proof of the above lemma essentially follows from that of \cite{Lou-Yong-2009}, Proposition 2.1. Based on the above, we further have the following result.

\bl{Lemma 4.9.} \sl Let $\g_i\in [0,1]$, $a_i(\cd)\in \cl{\sA[\g_i,\g_i]}^{\,H}$ $(i=1,2)$. Then for any $\th\in[0,1]$ and $e\in S^{n-1}$,
\bel{E111}b(\cd)\equiv a_1(\cd)-\th[a_1(\cd)-a_2(\cd)]-{\th(1-\th)[a_1(\cd)-a_2(\cd)]ee^\RT
[a_1(\cd)-a_2(\cd)]\over e^\RT\big\{a_2(\cd)+\th[a_1(\cd)-a_2(\cd)]\big\}e}\in\cl{\sA[\g,\g]}^{\,H},\ee
where $\g= (1-\th)\g_1+\th\g_2$. Consequently,
\bel{H(A)=A}\cH\[\,\cl{\sA[\a,\b]}^{\,H},\cl{\sA[\a,\b]}^{\,H};[0,1],S^{n-1}\]\subseteq
\cl{\sA[\a,\b]}^{\,H}.\ee
\el

\it Proof. \rm We need only to consider the case $\th\in(0,1)$, which is fixed below. Take an $e\in S^{n-1}$. For notation simplicity, denote $\th_1=1-\th$, $\th_2=\th$. For any $k\geq 1$,
denote
\bel{Qk}Q^1_k=\set{x\in\O\bigm|\big\{{k\lan x,e\ran}\big\}\in[\th,1)}, \q Q^2_k=\set{x\in\O\bigm|\big\{{k\lan x,e\ran}\big\}\in[0,\th)}
\ee
and define
\bel{b_k}
b_k(\cd)=\chi_{Q^1_k}(\cd)a_1(\cd)+\chi_{Q^2_k}(\cd)a_2(\cd).
\ee
Then
\bel{Qkw}
\chi_{Q^i_k}(\cd)\to \th_i\chi_\O(\cd),\qq \wi{L^2(\O)}, \q i=1,2.\ee
Moreover, by Lemma \ref{Lemma 4.8.}, $b_k(\cd)\hcon b(\cd)$.

\ms

On the other hand, for $i=1,2$, there is a sequence $\{E^i_j\}_{j\ges1}$ of measurable subsets of $\O$ such that
$\ds |E^i_j|=\g_i |\O|$,
\bel{E453}
\chi_{E^i_j}(\cd)\to \si_i(\cd), \q \wi{L^2(\O)},
\ee
and
\bel{E454}
a_j^i(\cd)=A_0+\chi_{_{E^i_j}}(\cd)(A_1-A_0)\hcon a_i(\cd).
\ee
By \eqref{Qkw},
\bel{gO}\sum^2_{i=1}\int_{Q^i_k}\si_i(x)dx\to\sum^2_{i=1}\int_\O\th_i\si_i(x) dx=\g|\O|,\ee
Thus, for $m\ges1$, we have $k_m\ges1$ such that (c.f. \eqref{rho})
\bel{rho-m}\rho\big(b_{k_m}(\cd),b(\cd)\big)\les{1\over m}\ee
and
\bel{si-m}\Big|\sum^2_{i=1}\int_{Q^i_{k_m}}\si_i(x)dx-\g|\O|\Big|\les{1\over m}.\ee
By Lemma \ref{Lemma 4.7N.}, as $j\to\infty$,
\bel{si-m}\sum^2_{i=1}\chi_{Q^i_{k_m}}(\cd)a^i_j(\cd)\hcon \sum^2_{i=1}\chi_{Q^i_{k_m}}(\cd)a_i(\cd)=b_{k_m}(\cd).\ee
Then, by \eqref{si-m} and \eqref{E453}, we have $j_m\ges1$ such that
\bel{rho-jm}\rho\big(\sum^2_{i=1} \chi_{Q^i_{k_m}}(\cd)a^i_{j_m}(\cd), b_{k_m}(\cd)\big)\les{1\over m}\ee
and
\bel{si-jm}\sum^2_{i=1}\Big|\int_{Q^i_{k_m}}\big(\chi_{E^i_{j_m}}(x)-\si_i(x)\big) dx\Big|\les{1\over m}.\ee
Denote
$$B_m(\cd)=\sum^2_{i=1} \chi_{Q^i_{k_m}}(\cd)a^i_{j_m}(\cd).$$
Then
$$B_m(\cd)=A_0 +\chi_{\O_m}(\cd)(A_1-A_0),\qq\O_m=\bigcup^2_{i=1}\big(Q^i_{k_m}\cap E^i_{j_m}\big).$$
By \eqref{rho-m} and \eqref{rho-jm}, $B_m(\cd)\hcon b(\cd)$. By \eqref{si-m} and \eqref{si-jm}, $\O_m\to \g|\O|$. Thus, it follows from
Lemma \ref{Lemma 4.7N.}{\rm(i)} that $b(\cd)\in \cl{\sA[\g,\g]}^{\,H}$,
proving our claim. \endpf

We further extend the above result to the following (replacing $e$ by $\xi(\cd)$).

\bl{Lemma 4.10.} \sl Assume $a(\cd),b(\cd)\in\cl{\sA[\a,\b]}^{\,H}$. Let $\O_1,\O_2,\cds,\O_m$ be a partition of $\O$ and
$$\xi(x)=\sum^m_{i=1}\chi_{_{\O_i}}(x)\xi_i,$$
with $\xi_i\in S^{n-1}$ ($1\les i\les m$). Then for any $\th\in [0,1]$,
\bel{EN610}
a(\cd)-\th[a(\cd)-b(\cd)]-{\th(1-\th)(a(\cd)-b(\cd))\xi(\cd)\xi(\cd)^\RT
(a(\cd)-b(\cd))\over\xi(\cd)^\RT\Big(b(\cd)+\th[a(\cd)-b(\cd)]\Big)
\xi(\cd)}\in\cl{\sA[\a,\b]}^{\,H}.\ee

\el

\it Proof. \rm For $m=1,2,\ldots,m$, denote
$$b_i(\cd)= a(\cd)-\th[a(\cd)-b(\cd)]-{\th(1-\th)[a(\cd)-b(\cd)]\xi_i\xi_i^\RT
[a(\cd)-b(\cd)]\over \xi_i^\RT\big\{b(\cd)+\th[a(\cd)-b(\cd)]\big\}\xi_i}.
$$
Since $a(\cd),b(\cd)\in\cl{\sA[\a,\b]}^{\,H}$, by Lemma \ref{Lemma 4.7.},
for $j=1,2$, there exist $\g_{j1},\g_{j2},\ldots,\g_{jm}\in [0,1]$  such that
$\ds a_j(\cd)\big|_{\O_i}\in\cl{\sA_{\O_i}[\g_{ji},\g_{ji}]}^{\,H}$, and
\bel{}\a|\O|\les\sum^m_{i=1}\g_{ji}|\O_i|\les\b|\O|,\ee
where we denote $a_1(\cd)=a(\cd)$, $a_2(\cd)=b(\cd)$ for notation simplicity.
By Lemma \ref{Lemma 4.9.}, $\ds
b_i(\cd)\big|_{\O_i}\in\cl{\sA_{\O_i}[\g_i,\g_i]}^{\,H}$ with $\g_i=(1-\th)\g_{1i}+\th \g_{2i}$ $(i=1,2,\ldots,m)$.
Since
5
$$\a|\O|\les\sum^m_{i=1}\g_i|\O_i|=(1-\th)\sum^m_{i=1}\g_{1i}|\O_i|+\th \sum^m_{i=1}\g_{2i}|\O_i|\les\b|\O|,$$
we get from Lemma \ref{Lemma 4.7.} that
$$\sum^m_{i=1}\chi_{\O_i}(\cd)b_i(\cd)\in\cl{\sA[\a,\b]}^{\,H}.$$
That is, \eqref{EN610} holds.  \endpf

Now, we are ready to prove Theorem \ref{T601}.

\ms

\it Proof of Theorem \ref{T601}. \rm Fix $b(\cd)\in\cl{\sA[\a,\b]}^{\,H}$. Let $\O_1,\O_2,\ldots,\O_m$ be a partition of $\O$ and
\bel{xi(x)}\xi(x)=\sum^m_{k=1}\chi_{_{\O_k}}(x)\xi_k,\ee
with $\xi_k\in S^{n-1}$ ($1\les k\les m$). Then, by Lemma \ref{Lemma 4.9.}, for any $\th\in (0,1)$,
\bel{EN610B}a_\th(\cd)\equiv\underline a(\cd)-\th[\underline a(\cd)-b(\cd)]-{\th(1-\th)[\underline a(\cd)-b(\cd)]
\xi(\cd)\xi(\cd)^\RT[\underline a(\cd)-b(\cd)]\over\xi(\cd)^\RT\Big(b(\cd)+\th[\underline a(\cd)-b(\cd)]\Big)\xi(\cd)}\in\cl{\sA[\a,\b]}^{\,H}.\ee
By the optimality of $\underline a(\cd)$, we have
$$\ba{ll}
\ns\ds\int_\O\lan\underline a(x)\nabla\underline y(x),\nabla\underline y(x)\ran\,dx=\l_{\underline a(\cd)}\les\l_{a_\th(\cd)}\\
\ns\ds=\inf_{y(\cd)\in W^{1,2}_0(\O)\atop\|y(\cd)\|_2=1}\int_\O\lan a_\th(x)\nabla y(x),\nabla y(x)\ran
\,dx\les\int_\O\lan a_\th(x)\nabla\underline y(x),\nabla\underline y(x)\ran\,dx.\ea$$
Therefore, by the minimality of $\underline a(\cd)$, one has
$$\ba{ll}
\ns\ds0\les\lim_{\th\to0^+}\int_\O\big\langle{a_\th(x)-\underline a(x)\over\th}\nabla\underline y(x),
\nabla\underline y(x)\big\rangle\,dx\\
\ns\ds\q=-\int_\O\big\langle\(\underline a(x)-b(x)+{(\underline a(x)-b(x))\xi(x)\xi(x)^\RT
(\underline a(x)-b(x))\over\xi(x)^\RT b(x)\xi(\cd)}\Big)\nabla\underline y(x),\nabla\underline y(x)\big\rangle\,dx.\ea$$
That is,
$$\ba{ll}
\ns\ds\int_\O\ip{\big[b(x)-\underline a(x)\big]\nabla\underline y(x),\nabla\underline y(x)}\,dx
\ges\int_\O{\big|\xi(x)^\RT\big[\underline a(x)-b(x)\big]\nabla\underline y(x)\big|^2\over
|b(x)^{1\over2}\xi(x)|^2}dx\\
\ns\ds\ges\int_\O\ip{b(x)^{-{1\over2}}\big[\underline a(x)-b(x)\big]\nabla\underline y(x),{b(x)^{1\over 2}\xi(x)\over|b(x)^{1\over2}\xi(x)|}}^2\,dx.\ea$$
The above is true for any $\xi(\cd)$ of form \eqref{xi(x)}. Then by approximation, we obtain
$$\ba{ll}
\ns\ds\int_\O\lan\big[b(x)-\underline a(x)\big]\nabla\underline y(x),\nabla\underline y(x)\ran\,dx\ges\int_\O\big|b(x)^{-{1\over 2}}(\underline a(x)-b(x))\nabla\underline y(x)\big|^2\,dx\\
\ns\ds=\int_\O\lan\big[\underline a(x)-b(x)\big]b(x)^{-1}\big[\underline a(x)-b(x)\big]\nabla
\underline y(x),\nabla\underline y(x)\ran\,dx\\
\ns\ds=\int_\O\lan\big[\underline a(x)b(x)^{-1}\underline a(x)-2\underline a(x)+b(x))\big]\nabla \underline y(x),\nabla\underline y(x)\ran\,dx.\ea$$
Therefore,
\eqref{EN601} holds.
Consequently, we have \eqref{EN602}.

Now, we show that \eqref{EN602} also implies \eqref{EN601}. Suppose that  \eqref{EN602} holds. For any $b(\cd)\in\cl{\sA[\a,\b]}^{\,H}$, there is a sequence $\O_k\subseteq\O$ with
$$\a|\O|\les|\O_k|\les\b|\O|,\qq k\ges1,$$
such that as $k\to+\infty$,
\bel{EN613B}b_k(\cd)\equiv\chi_{_{\O_k^c}}(\cd) A_0+\chi_{_{\O_k}}(\cd)A_1\hcon b(\cd).\ee
We can assume that
$$\chi_{_{\O_k}}(\cd)\to g(\cd),\qq\hb{weakly in }L^2(\O).$$
Then
$$b_k(\cd)^{-1}=\chi_{_{\O_k^c}}(\cd)A_0^{-1}+\chi_{_{\O_k}}(\cd)A_1^{-1}\to[1-g(\cd)]A_0^{-1}+g(\cd)A_1^{-1},\qq \mbox{weakly in } L^2(\O).$$
From (\ref{E205}), we get that
\bel{EN614}b(\cd)^{-1}\les[1-g(\cd)]A_0^{-1}+g(\cd)A_1^{-1}.\ee
Therefore,
$$\ba{ll}
\ns\ds\int_\O\lan\underline a(x)b(x)^{-1}\underline a(x)\nabla\underline y(x),\nabla\underline y(x)\ran\,dx\les\int_\O\lan\underline a(x)\big[(1-g(x))A_0^{-1}+g(x)A_1^{-1}\big]\underline a(x)\nabla \underline y(x),\nabla\underline y(x)\ran\,dx\\
\ns\ds=\lim_{k\to\infty}\int_\O\lan\underline a(x)\big[\chi_{_{\O_k^c}}(x)A_0^{-1}+\chi_{_{\O_k}}(x) A_1^{-1}\big]\underline a(x)\nabla\underline y(x),\nabla\underline y(x)\ran\,dx\les\int_\O\lan\underline a(x)\nabla\underline y(x),\nabla\underline y(x)\ran\,dx,\ea$$
where the last inequality follows from \eqref{EN602}. Hence, we get \eqref{EN601}. Then it follows that
\bel{sup=max}\ba{ll}
\ns\ds\sup_{b(\cd)\in\sA[\a,\b]}\int_\O\lan\underline a(x)b(x)^{-1}\underline a(x)\nabla
\underline y(x),\nabla\underline y(x)\ran\,dx\\
\ns\ds=\max_{b(\cd)\in\cl{\sA[\a,\b]}^{\,H}}\int_\O\lan\underline a(x)b(x)^{-1}\underline a(x)\nabla
\underline y(x),\nabla\underline y(x)\ran\,dx=\int_\O\lan\underline a(x)\nabla\underline y(x),\nabla\underline y(x)\ran\,dx=\underline\l.\ea\ee
Hence, \eqref{sup**} holds.

\ms

Next, in the case that $\a=0$ and $\b=1$, for any sub-domain $\O_\e\subseteq\O$ with $|\O_\e|=\e$, and $b(\cd)\in\cl{\sA[\a,\b]}^{\,H}$, let (noting Lemma \ref{Lemma 4.7.} (ii))
$$b_\e(\cd)=\chi_{_{\O_\e}}(\cd)b(\cd)+\chi_{_{\O_\e^c}}(\cd)\underline a(\cd)\in\cl{\sA[\a,\b]}^{\,H}.$$
Taking such a $b_\e(\cd)$ in the above, we obtain
$$\int_{\O_\e}\lan\underline a(x)b(x)^{-1}\underline a(x)\nabla\underline y(x),\nabla\underline y(x)\ran\,dx\les\int_{\O_\e}\lan\underline a(x)\nabla\underline y(x),\nabla\underline y(x)\ran\,dx.$$
Then, using Lebesgue's density theorem, we obtain \eqref{in bar cA}.
In paricular, \eqref{in A} holds.  Moreover, similar to that
\eqref{EN601} and \eqref{EN602} are
equivalent, \eqref{in A} is equivalent to \eqref{in A^H}.
\endpf

\subsection{Optimality system.}

Let us now take a closer look at \eqref{EN601}/\eqref{sup=max}.  Note that any $b(\cd)\in\sA[\a,\b]$ has the following
form:
$$b(\cd)=\chi_{_{\O_1^c}}(\cd)A_0+\chi_{_{\O_1}}(\cd)A_1,$$
for some $\O_1\subseteq\O$, with $\a|\O|\les|\O_1|\les \b|\O|$.
Then there is a sequence of $\O_k\subseteq \O$ such that $\a|\O|\les |\O_k|\les\b|\O|$,
$$\chi_{_{\O_k^c}}(\cd)A_0+\chi_{_{\O_k}}(\cd)A_1\hcon \underline a(\cd)$$
and
$$\chi_{_{\O_j}}(\cd)\to \uGs(\cd), \qq \mbox{weakly in } L^2(\O;[0,1]).$$
Then $\uGs(\cd)\in\Si[\a,\b]$ (see \eqref{Si[a,b]} for the definition) and by property (iv) in Subsection 4.1,
$$\ula(\cd)^{-1}\les A_0^{-1}+\uGs(\cd)(A_1^{-1}-A_0^{-1}),\qq x\in \O.$$
Consequently, \eqref{sup=max} is equivalent to the following:
$$\ba{ll}
\ns\ds\int_\O\ip{\big(A_0^{-1}+\uGs(x) (A_1^{-1}-A_0^{-1})\big)\ula(x)\nabla\underline y(x),\ula(x)\nabla\underline y(x)}dx\ges\int_\O\lan\underline a(x)\nabla\underline y(x),\nabla\underline y(x)\ran dx\\
\ns\ds=\sup_{\a|\O|\les|\O_1|\les\b|\O|}\int_\O\lan\underline a(x)\big[A_0^{-1}\chi_{_{\O_1^c}}(x)+A_1^{-1}\chi_{_{\O_1}}(x)\big]
\underline a(x)\nabla\underline y(x),\nabla\underline y(x)\ran dx\\
%
%
\ns\ds=\sup_{\si(\cd)\in \Si[\a,\b]}\int_\O\ip{\big(A_0^{-1}+\si(x) (A_1^{-1}-A_0^{-1})\big)\ula(x)\nabla\underline y(x),\ula(x)\nabla\underline y(x)} dx
.\ea$$
Therefore,
$$\ba{ll}
\ns\ds\int_\O\uGs(x)\ip{ (A_1^{-1}-A_0^{-1}) \ula(x)\nabla\underline y(x),\ula(x)\nabla\underline y(x)} dx\\
\ns\ds=\sup_{\si(\cd)\in \Si[\a,\b]}\int_\O\si(x)\ip{ (A_1^{-1}-A_0^{-1}) \ula(x)\nabla\underline y(x),\ula(x)\nabla\underline y(x)} dx.\ea$$
By denoting
$$h(x)=\ip{(A_0^{-1}-A_1^{-1})\ula(x) \na \uly(x),\ula(x) \na \uly(x)},$$
one sees that the above becomes
\bel{suph(x)}\int_\O\uGs(x)h(x)dx=\inf_{\si(\cd)\in\Si[\a,\b]}\int_\O\si(x)h(s)dx.\ee
Thus, $\underline\si(\cd)\in\Si[\a,\b]$ solves a maximization problem. For this problem, we have the following proposition.

\bp{maximization} \sl Let $\underline\si(\cd)\in\Si[\a,\b]$ satisfy \eqref{suph(x)}. Then there are two constants $\Psi$ and $\mu_0$ such that
$$\m_0\les0,\qq\m_0^2+\Psi^2=1,$$
$$\Big(\int_\O\uGs(x)dx-t\Big)\Psi\les 0,\qq \a| \O|\les t\les\b|\O|.
$$
and
$$\big(\m_0h(x)+\Psi\big)\uGs(x)=\max_{0\les\th\les1}\big(\m_0h(x)+\Psi\big)\th, \qq\eqae x\in\O.$$

\ep

\it Proof. \rm For any $\e>0$, define
$$\ba{ll}
\ns\ds F_\e(\si(\cd))=\Big\{\[\(\int_\O\big[\si(x)-\underline\si(x)\big]h(x)dx+\e\)^+
\]^2+\phi\Big(\int_\O\si(x)ds\Big)\Big\}^{1\over2},\ms\\
\qq\qq\qq\forall
\si(\cd)\in\h\Si=\big\{\si:\O\to[0,1]\bigm|\si(\cd)\hb{ is measurable }\big\},
\ea$$
where
$$\phi(s)=\min_{\a|\O|\les t\les\b|\O|}|s-t|^2,\qq s\in\dbR.$$
We have  $\ds \big|\phi^\prime(s)\big|=2\sqrt{\phi(s)}$ and
$$(t-s)\phi^\prime(s)\les0,\qq\a|\O|\les t\les\b|\O|.$$
Clearly,
$$F_\e(\cd)>0,\qq\forall\si(\cd)\in\h\Si,\qq F_\e(\underline\si(\cd))=\e.$$
By Ekeland's variational principle (\cite{Li-Yong-1995}), there exists a $\si_\e(\cd)\in \h\Si$ such that
$$\ba{ll}
\ns\ds\|\si_\e(\cd)-\underline\si(\cd)\|_2\les\sqrt\e,\\
\ns\ds F_\e(\si(\cd))+\sqrt\e\|\si(\cd)-\si_\e(\cd)\|_2\ges F_\e(\si_\e(\cd)),\qq\forall\si(\cd)\in\h\Si.\ea$$
Hence, for any $\d\in(0,1)$ and $\si(\cd)\in\h\Si$, one has
$$\si_\e^\d(\cd)\equiv\si_\e(\cd)+\d[\si(\cd)-\si_\e(\cd)]\in\h\Si,$$
and thus
$$\ba{ll}
\ns\ds-\sqrt\e\,\|\si(\cd)\|_2\les {F_\e\big(\si_\e^\d(\cd)\big)-F_\e\big(\si_\e(\cd)\big)\over\d}
={F_\e\big(\si_\e^\d(\cd)\big)^2-F_\e\big(\si_\e(\cd)\big)^2
\over\d\big[ F_\e\big(\si_\e^\d(\cd)\big)+F_\e\big(\si_\e(\cd)\big)\big]}\\
\ns\ds=\1n{1\over \big[F_\e\big(\si_\e^\d(\cd)\big)\1n+\1n F_\e\big(\si_\e(\cd)\big)\big]\d}
\Big\{\[\(\int_\O\big[\si_\e^\d(x)-\underline\si(x)\big]h(x)dx\1n+\1n\e\)^+
\]^2\3n-\1n\[\(\int_\O\big[\si_\e(x)-\underline\si(x)\big]h(x)dx\1n+\1n\e\)^+\]^2\\
\ns\ds\qq\qq+\phi\(\int_\O\si_\e^\d(x)dx\)-\phi\(\int_\O\si_\e(x)dx\)\Big\}\\
\ns\ds\to{1\over F_\e(\si_\e(\cd))}\[\(\int_\O[\si_\e(x)-\underline\si(x)]h(x)dx+\e\)^+
\int_\O[\si(x)-\si_\e(x)]h(x)dx\\
\ns\ds\qq\qq + {1\over 2}\phi^\prime\Big(\int_\O\si_\e(x)dx\Big)\int_\O[\si(x)-\si_\e(x)]dx\\
\ns\ds\equiv -\m_0^\e\int_\O[\si(x)-\si_\e(x)]h(x)dx-\Psi^\e\int_\O[\si(x)-\si_\e(x)]dx,\ea$$
with
$$(\m^\e_0)^2+(\Psi^\e)^2=1,\qq\m^\e_0\les 0,$$
$$\Big(\int_\O\si_\e(x)dx-t\Big)\Psi^\e\les0,\qq\a|\O|\les t\les\b|\O|.$$
Then along a subsequence, still denoted it by itself, we may let
$$(\m^\e_0,\Psi^\e)\to(\m_0,\Psi),\qq\m_0^2+\Psi^2=1,\qq\m_0\les0,$$
$$\Big(\int_\O\underline\si(x)dx-t\Big)\Psi\les0,\qq\a|\O|\les t\les\b|\O|,$$
and
$$\int_\O\big(\si(x)-\underline\si(x)\big)\big(\m_0h(x)+\Psi\big)dx\les0.$$
Hence, a standard argument applies to get
$$\big(\m_0h(x)+\Psi\big)\uGs(x)=\max_{0\les\th\les1}\big(\m_0h(x)+\Psi\big)\th, \qq\eqae x\in\O.
$$
We obtain our conclusions. \endpf

\ms

Now, we use the above result to make some further analysis on the optimal control $\underline\si(\cd)$ of Problem {\rm($\underline\L^H[\a,\b]$)}.

\ms

If $\mu_0=0$, then $\Psi\ne 0$ and we have
$\ds \uGs(\cd)\equiv 0 $ or $\ds \uGs(\cd)\equiv 1$.   That is $\ula(x)\equiv A_0$ or  $\ula(x)\equiv A_1$.

\ms

If $\mu_0\ne 0$, then we can suppose $\m_0=-1$ without loss of generality. Thus
$$\uGs(x)=\left\{\2n\ba{ll}
\ds1,\qq\eqae x\in (h(\cd)<\Psi),\\
\ns\ds0,\qq\eqae x\in (h(\cd)>\Psi).\ea\right.$$
This implies that
$$\ds|(h(\cd)\les\Psi)|\ges\int_\O\uGs(x)dx\ges\a|\O|,$$
and
$$\ds|(h(\cd)<\Psi)|\les\int_\O\uGs(x)dx\les\b|\O|.$$
Moreover, (when $\a<\b$) we can see that
$$\ba{ll}
\ns\ds\int_\O\uGs(x)dx>\a|\O|\qq\Ra\qq\Psi\les0;\\
\ns\ds\int_\O\uGs(x)dx<\b|\O|\qq\Ra\qq\Psi\ges0;\\
\ns\ds\Psi>0\qq\Ra\qq\int_\O\uGs(x)dx=\a|\O|;\\
\ns\ds\Psi<0\qq\Ra\qq\int_\O\uGs(x)dx=\b|\O|.\ea$$
For any $x\in(h(\cd)=\Psi)$, though it is possible that $\uGs(x)$ be any value of $[0,1]$, there are still some information could be used to determine
$\uGs(x)$. For example, if $|(h(\cd)\ges\Psi)|=\a|\O|$, then $\uGs(x)=0$ $\eqae(h(\cd)\ges\Psi)$.

\ms

On the other hand, when $\a |\O|<\int_\O \uGs(x) dx<\b |\O|$ (it will be the case if $\a=0,\b=1$ and neither $\chi_\O(\cd)A_0$ nor $\chi_\O(\cd)A_1$ is optimal), it should hold that $\Psi=0$. At this moment, on the set $(h(\cd)= \Psi)\equiv(h(\cd)=0)$,
\begin{equation}\label{E459}
\ip{A_1^{-1} \ula(x)\na \uly(x),\ula(x)\na \uly(x)}=\ip{ A_0^{-1}\ula(x)\na \uly(x),\ula(x)\na \uly(x)}.
\end{equation}
Hence, when $h(x)=0$, one has
$$\ba{ll}
\ns\ds\ip{A_1 \na \uly(x), \na \uly(x)}\ne \ip{ A_0^{-1}A_1\na \uly(x),A_1\na \uly(x)}\qq\Ra\qq\uGs(x)\ne 1;\\ [2mm]
\ns\ds\ip{ A_0\na \uly(x),\na \uly(x)}=\ip{A_1^{-1} A_0\na \uly(x),A_0\na \uly(x)}\qq\Ra\qq\uGs(x)\ne 0.\ea$$

\subsection{Maximization problem.}

Similar to Problem {\rm($\underline\L^H[\a,\b]$)}, it is natural to pose the following $H$-relaxation of Problem ($\bar\L[\a,\b]$).

\ms

\bf Problem ($\cl\L^H[\a,\b]$). \rm Find an $\bar a(\cd)\in\cl{\sA[\a,\b]}^{\,H}$ such that
\bel{4.61}\l_{\bar a(\cd)}=\sup_{a(\cd)\in\cl{\sA[\a,\b]}^{\,H}}\l_{a(\cd)}\ee

Any $\bar a(\cd)$ satisfying \eqref{4.61} is called an optimal control of Problem ($\cl\L^H[\a,\b]$), which is also called an optimal $H$-relaxed control of Problem ($\bar\L[\a,\b]$). Due to the properties of $H$-convergence, it is not hard to see that there are results for Problem ($\cl\L^H[\a,\b]$) parallel to the minimization problem. Let us state them here.

\bt{} \sl {\rm(i)} Problem $(\cl\L^H[\a,\b])$ admits an optimal control $\bar a(\cd)\in\cl{\sA[\a,\b]}^{\,H}$.

\ms

{\rm(ii)} Let $\bar a(\cd)\in\cl{\sA[\a,\b]}^{\,H}$ be an optimal control of Problem {\rm($\cl\L^H[\a,\b]$)} with $\bar y(\cd)\in W^{1,2}_0(\O)$ being the corresponding normalized principal eigenfunction. Then
\bel{inf*}\ba{ll}
\ns\ds\bar\l=\int_\O|\bar a(x)^{1\over2}\nabla\bar y(x)|^2dx=\inf_{b(\cd)\in\cl{\sA[\a,\b]}^{\,H}}\int_\O|b(x)^{-{1\over2}}\bar a(x)\nabla\bar y(x)|^2dx\\
\ns\ds\qq\qq\qq\qq\qq\q=\inf_{b(\cd)\in\sA[\a,\b]}\int_\O|b(x)^{-{1\over2}}
\bar a(x)\nabla\bar y(x)|^2dx,\ea\ee
and
\bel{in bar cA*}\int_\O\ip{\big[\bar  a(x)-\bar  a(x)b(x)^{-1}\bar a(x)\big]\nabla\underline y(x),\nabla\underline y(x)}\,dx\les0,\qq\all b(\cd)\in\cl{\sA[\a,\b]}^{\,H}.\ee
In paricular,
\bel{in cA*}\int_\O\ip{\big[\bar  a(x)-\bar  a(x)b(x)^{-1}\bar  a(x)\big]\nabla\underline y(x),\nabla\underline y(x)}\,dx\les0,\qq\all b(\cd)\in\sA[\a,\b].\ee
When $\a=0$ and $\b=1$, it holds that:
\bel{in A^H*}\lan\bar a(x)\nabla\bar y(x),\nabla\bar y(x)\ran\les\lan B^{-1}\bar a(x)\nabla\bar y(x),\bar a(x)\nabla\bar y(x)\ran,\q\ae x\in\O,~B\in\cl{\{A_0,A_1\}}^{\,H}.\ee
In paricular,
\bel{in A*}\lan\bar a(x)\nabla\bar y(x),\nabla\bar y(x)\ran\les\lan A_i^{-1}\bar a(x)\nabla\bar y(x),\bar a(x)\nabla\bar y(x)\ran,\q\ae x\in\O,~i=0,1.\ee
\et

The proof is omitted here. Also, one could derive (at least formally) the optimality system for the problem similar to the minimization problem.

\section{A Two-Dimensional Example.}

In this section, we present a two-dimensional example of Problem {\rm($\underline\L^H[0,1]$)}; Namely, $\a=0$, $\b=1$ and according to \eqref{ab}, we should assume that neither $A_0\les A_1$ nor $A_0\ges A_1$ holds. Since both $A_0$ and $A_1$ are positive definite, making a change of variables if necessary, without loss of generality, we may assume that
\bel{E501}A_0=\ppmatrix{\m_0&0\cr0&\m_1}, \q A_1=I,\ee
with $0<\m_0<1<\m_1$. Recall that this example is comparable with Example 3.3. Let $\underline a(\cd)\in\cl{\sA[0,1]}^{\,H}$ be an optimal control of Problem ($\underline\L^H[0,1]$). Then the following holds:
\bel{E503}\lan\underline a(x)\nabla\underline y(x),\nabla\underline y(x)\ran\ges\lan B^{-1}\underline a(x) \nabla\underline y(x),\underline a(x)\nabla\underline y(x)\ran,\q B=A_0,I,\qq\ae~x\in\O.\ee
To determine an optimal control $\underline a(\cd)$, let us make an observation. For given $x\in\dbR^n$, if we denote $\xi=\nabla\underline y(x)$ and $\bar A=\underline a(x)$, then \eqref{E503} reads
\bel{E504}\lan\bar A\xi,\xi\ran\ges\lan B^{-1}\bar A\xi,\bar A\xi\ran,\qq\forall B\in\{A_0,I\}.\ee
Or, equivalently,
\bel{EPP}\lan\bar A\xi,\xi\ran\ges\lan A_0^{-1}\bar A\xi,\bar A\xi\ran,\qq
\lan\bar A\xi,\xi\ran\ges\lan\bar A\xi,\bar A\xi\ran.\ee
Hence, it is natural to consider the following problem.

\ms

\bf Problem (P$_\xi$). \rm For $\xi\in \dbR^2$, find an $\bar A\in\cl{\{A_0,I\}}^{\,H}$ such that \eqref{E504} holds.

\ms

Here, we recall that
$$\cl{\{A_0,I\}}^{\,H}=\Big\{A\in M[\m_0,\m_1]\bigm|\chi_{_{\O}}(\cd)A
\in\cl{L^\infty\big(\O;\{A_0,I\}\big)}^{\,H}\Big\}.$$
For the above problem, we have the following interesting proposition.

\bp{Pxi} \sl {\rm(i)} If $\bar A\in\cl{\{A_0,I\}}^{\,H}$ is a solution to Problem {\rm ($\BP_\xi$)}, then
\bel{E506}
\lan\bar A\xi,\xi\ran\ges\lan B^{-1}\bar A\xi,\bar A\xi\ran,\qq\all B\in\cl{\{A_0,I\}}^{\,H}.\ee

{\rm(ii)} If both $\bar A$ and $\bar B$ are solutions of Problem ($\BP_\xi$), then $\bar A\xi=\bar B\xi$.

\ep

\it Prof. \rm (i) For any $B\in\cl{\{A_0,I\}}^{\,H}$, we can find a sequence $\O_k\subseteq\O$ such that
$$\chi_{_{\O_k^c}}(\cd)A_0+\chi_{_{\O_k}}(\cd)I\hcon\chi_{_\O}(\cd)B,\q {|\O_k|\over |\O|}\to\th,$$
for some $\th\in[0,1]$. By Property (vi) in Subsection 4.1,
$$B^{-1}\les(1-\th)A_0^{-1}+\th I.$$
Thus, \eqref{E506} follows.

\ms

(ii) Note that
$$\lan\bar A\xi,\xi\ran\ges\lan\bar B^{-1}\bar A\xi,\bar A\xi\ran,\qq
\lan\bar B\xi,\xi\ran\ges\lan\bar A^{-1}\bar B\xi,\bar B\xi\ran,$$
we have
$$\ba{ll}
\ns\ds\lan(\bar A^{-1}+\bar B^{-1})(\bar A-\bar B)\xi,(\bar A-\bar B)\xi\ran=\lan\bar B^{-1}\bar A\xi,\bar A\xi\ran-\lan\bar A\xi,\xi\ran+\lan\bar A^{-1}\bar B\xi, \bar B\xi\ran-\lan\bar B\xi,\xi\ran\les0.\ea$$
Therefore, it must hold that $\bar A\xi=\bar B\xi$. \endpf

The above tells us that to meet the necessary conditions for optimal controls of Problem ($\underline\L^H[0,1]$) on $\cl{\sA[0,1]}^{\,H}$, it suffices to find, for almost each $x\in\O$, an $\bar A\in\cl{\{A_0,I\}}^{\,H}$ satisfying
\bel{E507}\lan\bar A\nabla\underline y(x),\nabla\underline y(x)\ran\ges\lan B^{-1}\bar A\nabla\underline y(x),\bar A\nabla\underline y(x)\ran,\q B=A_0,I.\ee
Part (ii) of above proposition means that although $\bar A$ might not be unique, $\bar A\xi$ is unique. Thus, if one can solve Problem ($\BP_\xi$) successfully for each $\xi\in \dbR^2$, then we obtain a map $\bar A:\dbR^2\to\dbS^2$. Then
\bel{optimal a}\underline a(x)=\bar A(\nabla\underline y(x)),\qq x\in\O\ee
gives an optimal control, where $\underline y(\cd)$ is a solution to the closed-loop system:
\bel{closed-loop}\left\{\2n\ba{ll}
\ds-\nabla\cd\(\bar A\big(\nabla\underline y(x)\big)\nabla\underline y(x)\)=\underline\l\,\underline y(x),\qq x\in\O,\\
\ns\ds\underline y\big|_{\pa\O}=0.\ea\right.\ee
We will see later it is a nonlinear eigenvalue problem.

\ms

Now, for given $\xi=(\xi_1,\xi_2)^\RT\ne 0$, we try to find a solution $\bar A$ of Problem ($\BP_\xi$).

\ms

Let us introduce the following partition of $\dbR^2$:
\begin{equation}\label{E508}\ba{l}
\ns\ds E_{A_0}=\Big\{(\xi_1,\xi_2)^\RT\ne 0\bigm|\xi_2^2\les{\m_0(1-\m_0)\over\m_1(\m_1-1)}\xi_1^2\Big\},\\
\ns\ds E_I=\Big\{(\xi_1,\xi_2)^\RT\ne 0\bigm|\xi_2^2\ges{\m_1(1-\m_0)\over\m_0(\m_1-1)}\xi_1^2\Big\},\\
\ns\ds E_+=\set{(\xi_1,\xi_2)^\RT\in\dbR^2\bigm|{\m_0(1-\m_0)\over\m_1(\m_1-1)}\xi_1^2<\xi_2^2<{\m_1(1-\m_0)
\over\m_0(\m_1-1)}\xi_1^2, \q \xi_1\xi_2>0},\\
\ns\ds E_-=\set{(\xi_1,\xi_2)^\RT\in\dbR^2\bigm|{\m_0(1-\m_0)\over\m_1(\m_1-1)}\xi_1^2<\xi_2^2<{\m_1(1-\m_0)
\over\m_0(\m_1-1)}\xi_1^2, \q \xi_1\xi_2<0}.\ea
\end{equation}
These sets are illustrated in the following figure. Clearly, $E_{A_0}, E_I, E_+, E_-$ are non-empty, mutually disjoint and $\ds \dbR^2=E_{A_0}\cup E_I\cup E_+\cup E_-\cup\set{0}$. Note that $E_{A_0}\cup\{0\}$ and $E_I\cup\{0\}$ are closed, and $E_\pm$ are open.

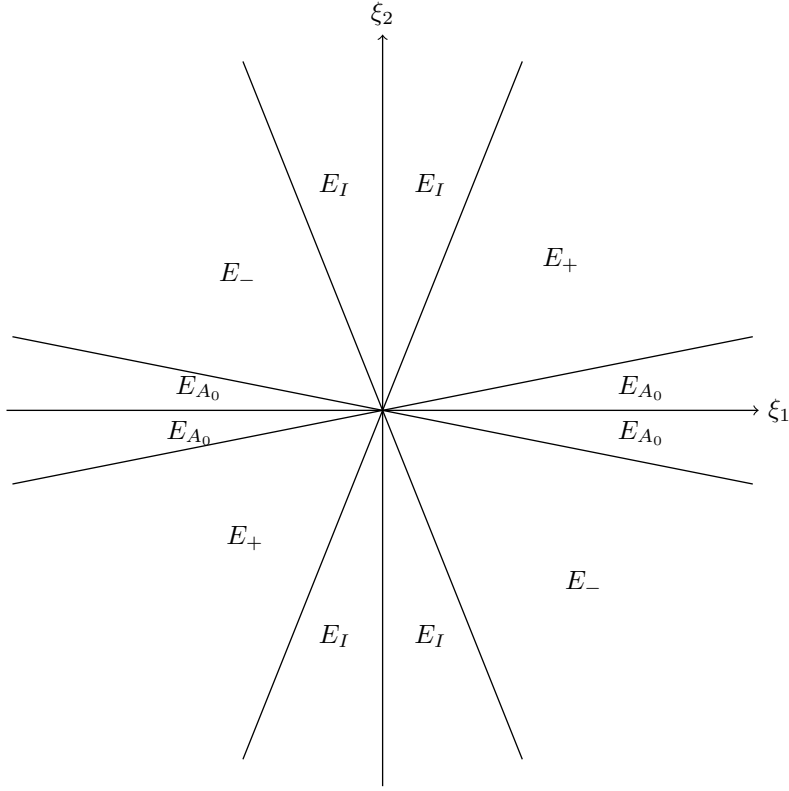
\begin{center}
\begin{tikzpicture}[line width=0.5pt]
\draw[->] (-5,0)--(5,0);
\draw[->] (0,-5)--(0,5);
\draw (-1.857,-4.6424)--(1.857,4.6424);
\draw (1.857,-4.6424)--(-1.857,4.6424);

\draw (-4.921,-0.98)--(4.921,0.98);
\draw (-4.921,0.98)--(4.921,-0.98);

\node[above] at (0,5) {$\xi_2$};
\node[right] at (5,0) {$\xi_1$};
\node[right] at (0.3,3) {$E_I$};
\node[left] at (-0.3,3) {$E_I$};
\node[right] at (0.3,-3) {$E_I$};
\node[left] at (-0.3,-3) {$E_I$};

\node[right] at (3,0.3) {$E_{A_0}$};
\node[right] at (3,-0.3) {$E_{A_0}$};
\node[left] at (-2,0.3) {$E_{A_0}$};
\node[right] at (-3,-0.3) {$E_{A_0}$};

\node[right] at (2,2) {$E_+$};
\node[right] at (-2.2,-1.7) {$E_+$};
\node[right] at (2.3,-2.3) {$E_-$};
\node[right] at (-2.3,1.8) {$E_-$};

\end{tikzpicture}
\end{center}

\no With a little calculation, one can see the following:

\ms

$\bullet$ $A_0$ is a solution of Problem ($\BP_\xi$)  if and only if $\ds \ip{A_0\xi,\xi}\ges\ip{A_0\xi,A_0\xi}$, i.e., $\xi\in E_{A_0}$.

\ms

$\bullet$ $I$ is a solution of Problem ($\BP_\xi$)  if and only of
$\ds \ip{\xi,\xi}\ges\ip{A_0^{-1}\xi,\xi}$, i.e., $\xi\in E_I$.

\ms

$\bullet$ Neither $A_0$ nor $I$ is a solution of Problem ($\BP_\xi$) if and only if $\xi\in E_+\cup E_-$, i.e.,
\begin{equation}\label{E509}
\ip{A_0\xi,\xi}<\ip{A_0\xi,A_0\xi}, \q \ip{\xi,\xi}< \ip{A_0^{-1}\xi,\xi}.
\end{equation}

\ms

Now, let $\ds\xi=(\xi_1,\xi_2)^\RT\in E_+\cup E_-$, by definition, $\xi_1,\xi_2\ne0$. Let $\eta=\bar A\xi$. Since $\bar A\in\cl{\{A_0,I\}}^{\,H}\setminus\{A_0,I\}$, by the proof of Proposition 5.1, (i), there exists a $\g\in(0,1)$ such that
\bel{E510}\bar A^{-1}\les(1-\g)A_0^{-1}+ \g I.\ee
Then
$$\lan\bar A^{-1}\eta,\eta\ran\les(1-\g)\lan A_0^{-1}\eta,\eta\ran+ \g \lan\eta,\eta\ran.$$
On the other hand, (\ref{EPP}) implies
$$\lan\bar A^{-1}\eta,\eta\ran\ges\lan A_0^{-1}\eta,\eta\ran, \q
\lan\bar A^{-1}\eta,\eta\ran\ges\lan\eta,\eta\ran.$$
Thus, it should hold that
\bel{E511}\lan\bar A^{-1}\eta,\eta\ran=\lan A_0^{-1}\eta,\eta\ran=\lan \eta,\eta\ran,\ee
which coincides with \eqref{E459}. From the second equality in \eqref{E511}, a direct calculation shows
\bel{E512}\eta=C_\xi \ppmatrix{\e\sqrt{1-s}\cr\sqrt s},\ee
with $C_\xi\ne0$ and $\e=\pm1$, where
\bel{E513}s={(1-\m_0)\m_1\over \m_1-\m_0}\in (0,1).\ee
Let us now determine $C_\xi$, $\e$ and $\g$. We have
$$C_\xi^2 =\ip{\eta,\eta}=\ip{\bar A^{-1}\eta, \eta}=\ip{\xi,\eta}=C_\xi \ip{ \ppmatrix{\e \sqrt{1-s}\cr \sqrt s},\xi}.$$
Thus
\bel{E514}
C_\xi=\e \xi_1\sqrt{1-s}+\xi_2\sqrt s\ee
and
\bel{E515}\ba{ll}
\ns\ds\bar A\xi=\eta=C_\xi\ppmatrix{\e \sqrt{1-s}\cr \sqrt s}= \ppmatrix{\e \sqrt{1-s}\cr \sqrt s}\ppmatrix{\e \sqrt{1-s}\cr \sqrt s}^\RT \xi=\ppmatrix{1-s & \e\sqrt{s(1-s)}\cr \e\sqrt {s(1-s)} & s} \xi\equiv G\xi.\ea\ee
On the other hand, by \eqref{E510} and \eqref{E511}, we have
$$\Big|\Big((1- \g)A_0^{-1}+ \g I-\bar A^{-1}\Big)^{1\over 2}\eta\Big|^2=\ip{\Big((1- \g)A_0^{-1}+ \g I-\bar A^{-1}\Big)\eta,\eta}=0.$$
This implies
\begin{equation}\label{E516}
\Big((1- \g)A_0^{-1}+ \g I-\bar A^{-1}\Big)\eta=0.
\end{equation}
Thus,
\begin{equation}\label{E517}
\xi= \bar A^{-1}\eta=\Big((1-\g) A_0^{-1}+\g I\Big)\eta= C_\xi  \ppmatrix{\e\Big({1-\g\over \m_0}+ \g \Big)\sqrt{1-s}\cr \Big({1-\g\over \m_1}+\g \Big)\sqrt s}.
\end{equation}
Hence, $\sgn(C_\xi)=\sgn(\xi_2)$, which can be obtained from \eqref{E514} and $x\in E_+\cup E_-$ too. Consequently,
\begin{equation}\label{E518}
 \qq \e=\sgn(\xi_1\xi_2).
\end{equation}
Moreover, it follows from $\ds \xi_2=C_\xi \Big({1-\g\over \m_1}+\g \Big)\sqrt s$ that
\begin{equation}\label{E519}
 \qq \g={|\xi_2|\sqrt{\m_1\over 1-\m_0}-|\xi_1|\sqrt{\m_0\over \m_1-1}\over |\xi_1|\sqrt{\m_0(\m_1-1)}+|\xi_2|\sqrt{\m_1(1-\m_0)}}.
\end{equation}
One can verify that the above $\g$ belonging to $(0,1)$ is equivalent to $\xi\in E_+\cup E_-$. Actually, it is a one-to-one mapping from $\ds \set{(\xi_1,\xi_2)^\RT\in S^1\cap(E_+\cup E_-)\big| \xi_1>0, \xi_2>0}$ to $(0,1)$.

\ms

Next, by (ii) of Proposition \ref{Pxi}, $\eta\equiv \bar A \xi$ only depends on $\xi$ (independent of the solution $\bar A$). We now find an $\bar A\in\cl{\{A_0,I\}}^{\,H}$ such that $\bar A\xi=\eta$. To this end, we try to find an solution in $\G(A_0,I)$. By  careful calculation, we find such a solution
$\bar A$ as the following:
\bel{E520}\bar A=I+(1-\g) Q^{-1},\qq Q=(A_0-I)^{-1}+\g (I-G).\ee
Denote $H=I-G$. Then $H\in\dbS^n$ with $H\ges 0$ and $\tr(H)=1$. Thus $\bar A\in \G(A_0,I)$. Moreover,
$H^2=H$. Let us verify $\bar A$ defined by \eqref{E520}  really satisfies $\bar A \xi=\eta$.
We have
$$
(A_0-I)QH=H+\g (A_0-I)H.
$$
Thus
$$\ba{ll}
\ns\ds A_0^{-1}(A_0-I)Q\big(\bar A\xi-\eta\big)=A_0^{-1}(A_0-I)Q\Big(H+(1-\g)Q^{-1}\Big)\xi\\
\ns\ds=A_0^{-1}\Big(H+\g (A_0-I)H+(1-\g)(A_0-I)\Big)\xi\\
\ns\ds=A_0^{-1}\Big(A_0-(1-\g)G-\g A_0 G\Big)\xi=\xi-(1-\g)A_0^{-1}\eta-\g\eta=0.\ea$$
Thus, along with  $G, \g, \e$ being given by \eqref{E515}, \eqref{E518} and \eqref{E519}, (when $\xi\in E_+\cup E_-$) we could get $\bar A$ by \eqref{E520}. Generally, we can choose
\bel{E538}\bar A(\xi)=\left\{\begin{array}{ll}
\ds  I,& \xi =0,\\
\ds A_0, &  \xi  \in E_{A_0},\\
\ds I,& \xi  \in E_I,\\
 \ds I+(1-\g)\Big((A_0-I)^{-1}+\g (I-G)\Big)^{-1},& \xi \in E_+\cup E_-.
\end{array}
\right.
\end{equation}
Therefore \eqref{closed-loop} is a nonlinear eigenvalue problem. We now simplify it. Noting \eqref{E515}, we see that
\bel{E529}
\ula(x)\na \uly(x)=\bar A(\na \uly(x))\na \uly(x)=\BF(\na \uly(x)),\qq\ae\,x\in\O,
\ee
where $\BF:\dbR^n\to \dbR^n$ is defined as
\begin{equation}
\label{E530}
\BF(\xi)=\left\{\begin{array}{ll}
\ds 0, & \xi=0,\\
\ds A_0\xi, &  \xi  \in E_{A_0},\\
\ds\xi,& \xi  \in E_I,\\
 \ds G_{_+}\xi,& \xi \in E_+,\\
\ns\ds G_{_-}\xi,& \xi \in E_-,
\end{array}
\right.
\end{equation}
with
\begin{equation}\label{E531}
G_{_\pm}= \ppmatrix{1-s & \pm\sqrt{s(1-s)}\cr \pm \sqrt {s(1-s)} & s}.
\end{equation}
Since $\BF(\na\uly(x))=\ula (x)\na \uly(x)$, $\uly(\cd)$ solves
\bel{E536}\left\{\2n\ba{ll}
\ds-\nabla\cd\big(\BF(\na \uly(x))\big)=\underline\l\underline y(x),\qq\eqin\O,\\
\ns\ds\underline y\big|_{\pa\O}=0.\ea\right.\ee
\ms
Although $G_\pm$ are singular (0 is an eigenvalue of $G_\pm$), noting that $\BF(\xi)=\bar A(\xi)\xi$, we still have
\begin{equation}
\label{E532}
\m_0|\xi|^2\les\ip{\BF(\xi),\xi}\les\m_1|\xi|^2, \qq\all \xi\in \dbR^n.
\end{equation}
On the other hand, it is not difficult to verify that $\ds\ip{\BF(\cd),\cd}$ is convex in $\dbR^n$ and
\begin{equation}\label{E535}
\big|\BF(\xi)-\BF(\tilde \xi)\big|\les\m_1 |\xi-\tilde \xi|, \qq\all \xi,\tilde \xi\in \dbR^n.
\end{equation}
Consequently,
 there is $\tilde y(\cd)\in W^{1,2}_0(\O)$ such that
$$
 \tilde \l\equiv {\int_\O\ip{\BF(\na \tilde y(x)),\na \tilde y(x)} \, dx\over \int_\O\|\tilde y(x)|^2\, dx}=\inf_{y(\cd)\in W^{1,2}_0(\O)\atop y(\cd)\ne 0}{\int_\O\ip{\BF(\na y(x)),\na y(x)} \, dx\over \int_\O\|y(x)|^2\, dx}.
$$
We have
$$
\left\{\2n\ba{ll}
\ds-\nabla\cd\big(\BF(\na \tilde y(x))\big)=\tilde\l\tilde y(x),\qq\eqin\O,\\
\ns\ds\tilde y\big|_{\pa\O}=0.\ea\right.
$$
Moreover, let $\ds\tilde a(\cd)=\bar A(\tilde y(\cd))$. Then $\tilde a(\cd)\in \cl{\sA[0,1]}^{\,H}$ and
$$\left\{\2n\ba{ll}
\ds-\nabla\cd\big(\tilde a(x)\nabla \tilde y(x)\big)=\tilde \l \tilde y(x),\qq\hb{in }\O,\\
\ns\ds \tilde y\big|_{\pa\O}=0,\ea\right.
$$
which implies $\underline \l\les\l_{\tilde a(\cd)}\les\tilde\l\les\underline \l$. This implies
\begin{equation}\label{E539}
\underline\l= {\int_\O\ip{\BF(\na \uly(x)),\na \uly(x)} \, dx\over \int_\O\|\uly(x)|^2\, dx}=\inf_{y(\cd)\in W^{1,2}_0(\O)\atop y(\cd)\ne 0}{\int_\O\ip{\BF(\na y(x)),\na y(x)} \, dx\over \int_\O\|y(x)|^2\, dx}.
\end{equation}
Therefore,  $\uly(\cd)$ is an optimal state of  Problem {\rm($\underline\L^H[0,1]$)} if and only if it is a minimizer of
$$
\sF(y(\cd))={\int_\O\ip{\BF(\na y(x)),\na y(x)} \, dx\over \int_\O\|y(x)|^2\, dx}
$$
over $W^{1,2}_0(\O)\setminus \set{0}$.

\ms

The results of this section can be summarized as follow: To get a solution of Problem {\rm($\underline\L^H[0,1]$)}, one can first find a nontrivial solution of \eqref{E536} with the smallest positive number $\underline\l$, or
equivalently, find a minimizer of $\sF(\cd)$ over $W^{1,2}_0(\O)\setminus \set{0}$. Then define $\ula(x)=\bar A(\na\uly(x))$, getting a solution of Problem {\rm($\underline\L^H[0,1]$)}.

\section{Concluding Remarks.}

We have investigated the maximization and minimization problems of the principal eigenvalue of elliptic operators with the Dirichlet boundary condition. The control appears in the diffusion matrix (the leading coefficient). These problems are well-motivated by composite material design to optimize the heat conduct property of the material (cooling down as quick as possible, or preserving the temperature as long as possible). For maximization problem, due to the concavity of the principle eigenvalue as a functional of the leading coefficient, as long as the control set is convex, optimal control and its characterization can be obtained easily. When the control set is not convex, we introduce the usual convexification to guarantee the existence of an optimal relaxed control. Then some necessary conditions can also be obtained. From an example, we see that uniformly mixing two material might not be optimal in the maximization problem.

\ms

For minimization problem, the situation is much more complicated due to the concavity of the principle eigenvalue as a functional of the control. We adopt the $H$-convergence so that the existence of the $H$-relaxed optimal control could be guaranteed. Instead of looking at the most general situation, we concentrate on the case of the lamination of two material whose diffusibility matrices are given. Some interesting necessary conditions are derived. It is worthy of pointing out that even both two material have their diagonal diffusibility matrices, the optimal diffusibility matrix could be non-diagonal. Such a situation has been exhibited through an illustrative example in Section 5.

\end{document}